\documentclass[article]{amsart}
\usepackage{cases}
\usepackage{caption}
\usepackage{amsmath, amsfonts, pifont, amssymb}
\usepackage{verbatim}
\usepackage[bookmarks=true]{hyperref}
\usepackage{tikz}
\usepackage{tkz-euclide}
\usetikzlibrary{positioning}
\usetikzlibrary{calc}
\usetikzlibrary{snakes,arrows}
\usetikzlibrary{quotes,angles}
\usepackage{geometry}
\usepackage{color}

\newtheorem{theorem}{Theorem}[section]
\newcommand{\rt}{\R\times\T_{\lambda}}
\newtheorem{lemma}[theorem]{Lemma}
\newtheorem{proposition}[theorem]{Proposition}
\newtheorem{corollary}[theorem]{Corollary}

\newcommand{\R}{\mathbb{R}}

\newcommand{\rz}{\R\times\Z_{1/\lambda}}
\newcommand{\beq}{\begin{equation}}
\newcommand{\eeq}{\end{equation}}
\newcommand{\beqq}{\begin{equation*}}
\newcommand{\eeqq}{\end{equation*}}

\theoremstyle{definition}
\newtheorem{definition}[theorem]{Definition}

\newcommand{\RT}{\R\times\T}
\theoremstyle{remark}
\newtheorem{remark}[theorem]{Remark}

\numberwithin{equation}{section}
\newcommand{\T}{\mathbb{T}}
\newcommand{\Z}{\mathbb{Z}}

\begin{document}
\title[cubic NLS on $\R\times\T$]{Improved global well-posedness for the cubic NLS on two-dimensional waveguide $\R\times\T$}

\author{Qionglei Chen}
\address{Qionglei Chen
\newline \indent Institute of Applied Physics and Computational Mathematics,
Beijing 100088,\ P. R. China,
\newline\indent
National Key Laboratory of Computational Physics, Beijing 100088, China}
\email{Chen\_qionglei@iapcm.ac.cn, }	
\author{Yilin Song}
\address{Yilin Song
\newline \indent The Graduate School of China Academy of Engineering Physics,
Beijing 100088,\ P. R. China}
\email{songyilin21@gscaep.ac.cn}
\author{Kailong Yang}
	\address{Kailong Yang
		\newline \indent National Center for Applied Mathematics in Chongqing, Chongqing Normal University, Chongqing,  \\
		\indent China. }
	\email{ykailong@mail.ustc.edu.cn}
\author{Ruixiao Zhang}
\address{Ruixiao Zhang
\newline \indent The Graduate School of China Academy of Engineering Physics,
Beijing 100088,\ P. R. China}
\email{zhangruixiao21@gscaep.ac.cn}
	
\author{Jiqiang Zheng}
\address{Jiqiang Zheng
\newline \indent Institute of Applied Physics and Computational Mathematics,
Beijing, 100088, China.
\newline\indent
National Key Laboratory of Computational Physics, Beijing 100088, China}
\email{zheng\_jiqiang@iapcm.ac.cn, zhengjiqiang@gmail.com}	

\begin{abstract}
In this article, we show that the solution to defocusing cubic nonlinear Schr\"odinger equation (NLS) posed on the two-dimensional waveguide 
\begin{align*}
i\partial_tu+\Delta_{\R\times\T}u=|u|^2u
\end{align*}
is globally well-posed in $H^s(\R\times\T)$ with $s>\frac{1}{2}$. The proof is based on the $I$-method. Inspired by Colliander-Keel-Staffilani-Takaoka-Tao [Discrete Contin. Dyn. Syst. 21 (2008), 665-686], we  construct  the modified energy to  improve the energy increment. The main difficulty lies in controlling the resonant interactions caused by the modified energy. To this end, we establish refined bilinear Strichartz estimates with angular truncation on the rescaled waveguide, thereby generalizing results previously obtained by Takaoka [J. Differ. Equa. 394 (2024), 296-319]. Furthermore, we demonstrate polynomial growth of  $H^s$ with $\frac{1}{2} < s < 1$. Our result extends the recent work of Deng-Fan-Yang-Zhao-Zheng [J. Func. Anal. 287 (2024), 110595].  
\end{abstract}	
\maketitle
\section{Introduction}
In this paper, we investigate the defocusing cubic NLS posed on the two-dimensional waveguide manifold $\R\times\T$
\begin{align}\label{NLS-intro}
\begin{cases}
i\partial_tu+\Delta_{\R\times\T}u=|u|^2u,\quad(t,x,y)\in\R_{t}\times\R\times\T,\\
u(0,x,y)=u_0(x,y)\in H^s(\RT)
\end{cases}
\end{align}
where $u:\R_t\times\R\times\T\to \mathbb{C}$ is the unknown function. The waveguide manifold $\R\times\T$ is known as the semi-periodic space. The Cauchy problem \eqref{NLS-intro} has the Hamiltonian structure and its solution  enjoys the following two conservation laws
\begin{align*}
\mbox{Mass: }&M(u(t))=\int_{\R\times\T}|u(t,x,y)|^2dxdy=M(u_0),\\
\mbox{Energy: }&E(u(t))=\int_{\R\times\T}\left(\frac{1}{2}|\nabla u(t,x,y)|^2+\frac{1}{4}|u(t,x,y)|^4\right)dxdy=E(u_0).
\end{align*}
	
The local well-posedness of $\eqref{NLS-intro}$ is influenced by the choice of the spatial domain. In the case of $\R^2$, this equation exhibits invariance under the following scaling transformation, which in turn gives rise to the concept of criticality for the Cauchy problem associated with \eqref{NLS-intro}: 
\begin{align*}
u(t,x)\rightarrow  u_\lambda(t,x)=\lambda u(\lambda^2t,\lambda x),\quad\forall\lambda>0.
\end{align*} 
Indeed, one can verify that the $L^2$ norm remains invariant under the above scaling transformation.
Strichartz estimates enable us to derive the local existence of $\eqref{NLS-intro}$ in $H^s(\Bbb{R}^2)$ for $s\geqslant 0$.  For further details regarding the local theory of NLS on \( \mathbb{R}^2\), please refer to Cazenave\cite{Caze}.
	
The fundamental work on local well-posedness for the NLS on the torus was initiated by Bourgain\cite{Bour1993}. In this pioneered paper, he first established  periodic Strichartz estimates for the linear Schr\"{o}dinger equation through analytical number theory. These Strichartz estimates are  optimal on $\T^2$. Subsequently,  the local well-posedness of \eqref{NLS-intro} in $H^s(\T^2)$ for $s > 0$ was derived using these estimates. For the semi-periodic case, i.e. $\R\times\T$, Tzvetkov-Takaoka\cite{TT} demonstrated the $L_{t,z}^4$ Strichartz estimate for frequency-localized initial data without any loss of derivatives. Using this Strichartz estimate, they established the local well-posedness of \eqref{NLS-intro} in $H^s(\R\times\T)$ with $s\geq0$.

According to the local well-posedness theory and conservation of mass and energy, the global well-posedness of \eqref{NLS-intro} in energy space $H^1(\RT)$ is quite direct. Therefore, the global well-posedness in $H^s$ with $s<1$ is more interesting since the equation is $L^2$ invariant in the Euclidean space. We also expect the critical regularity is $L^2$ in the semi-periodic setting. For the case of $\R^2$, Bourgain first developed the Fourier truncation method to prove the $\eqref{NLS-intro}$ is globally well-posed in $H^s(\R^2)$ with $s>\frac23$. Colliander-Keel-Staffilani-Takaoka-Tao\cite{CKSTT2002Almost}  initially exploited the I-method to obtain the global well-posedness for $\eqref{NLS-intro}$ in $H^s(\R^2) $ with $s>\frac47$. By taking advantage of  the resonant decomposition and $I$-method, they further improved the result of \cite{CKSTT2002Almost}  to $s>\frac12$ in \cite{CKSTT2008Res}. In 2016, Dodson\cite{Dodson} established the long-time Strichartz estimates and employed the concentration-compactness argument developed by Kenig-Merle \cite{KM} to prove the global well-posedness and scattering in $L^2(\R^2)$. For the case of $\T^2$,  De Silva-Pavlovi\'{c}-Staffilani-Tzirakis\cite{Silva2007Global} proved the global well-posedness of the cubic NLS in $H^s(\Bbb{T}^2)$ with $s>\frac{2}{3}$ by using I-method in conjunction with refined bilinear Strichartz estimates on rescaled tori $\T_{\lambda}^2=(\lambda\T)^2$,
\begin{align}\label{bil-1}
\bigg\Vert \eta(t)U_\lambda\phi_1\cdot\eta(t)U_\lambda\phi_2\bigg\Vert_{L_{t,x}^2(\R\times\T_\lambda^2)}\lesssim N_2^\varepsilon\left(\frac{1}{\lambda}+\frac{N_2}{N_1}\right)^\frac{1}{2}\Vert\phi_1\Vert_{L^2(\T_\lambda^2)}\Vert\phi_2\Vert_{L^2(\T_\lambda^2)},
\end{align} 
where $1\leqslant N_2\ll N_1$, $\lambda\geqslant1$, $\eta(t)\in C_0^\infty([-4,4])$, $\phi_j\in L^2(\lambda\Bbb{T}^2)$ is spectrally localized at $\left|k_j\right|\sim \lambda N_j$($j=1,2$) and 
\begin{align*}
U_\lambda(t)u_0(x)=\frac{1}{\lambda}\sum_{k\in\frac{1}{\lambda}\Bbb{Z}^2}e^{2\pi ik\cdot x-4\pi it\left|k\right|^2}\hat{u}_0(k).
\end{align*} 
We note that when $\lambda\rightarrow\infty$, this bilinear estimates $\eqref{bil-1}$ is consistent with the bilinear Strichartz estimates in $\Bbb{R}^2$ up to the $\varepsilon$ loss. Schippa\cite{Schippa} utilized the $I$-method along with the refined Strichartz estimate developed in \cite{Refined} to show that the cubic NLS posed on $\T^2$ is globally well-posed in $H^s(\T^2)$ for $s>\frac35$. Herr-Kwak\cite{Herr,Herr2} established the global existence for solution to $\eqref{NLS-intro}$ in $H^s(\T^2)$ for $s>0$. For  NLS on general waveguide manifolds  $\R^d\times\T^m$, we  refer to  \cite{Luo-JMPA,Luo-MA,Zhao-JDE,Zhao-Zheng}. 

Recently, Deng-Fan-Yan-Zhao-Zheng\cite{DFYZZ} established a sharp local-in-time bilinear Strichartz estimate without any loss of derivatives and then combine with $I$-method to show the global well-posedness for solution to \eqref{NLS-intro} in $H^s(\RT)$ for $s>\frac35$. The main aim of this paper is to further improve the result in \cite{DFYZZ} by making use of the resonant decomposition technique developed in \cite{CKSTT2008Res}.
\begin{theorem}\label{thm1}
The solution to \eqref{NLS-intro} is globally well-posed in $H^s(\RT)$ with $s>\frac{1}{2}$. Moreover, we also have the polynomial growth of Sobolev norms 
\begin{equation*}
\|u(T)\|_{H^s(\RT)} \lesssim T^{\frac{s(1-s)}{4s-2}+},\quad \forall T > 1.
\end{equation*}
\end{theorem}
\begin{remark}
The polynomial growth of the higher-order Sobolev norms is another interesting problem, which is tightly linked to the weak turbulence theory.		Recently, Takaoka \cite{Takaoka} applied the upside-down I-method as well as the resonant decomposition to obtain the upper bound of the growth of higher order Sobolev norms with $s>1$:
\begin{equation*}
\| u(T) \|_{H^s(\RT)} \lesssim (1+T)^{\frac{s-1}{2}+},\quad \forall T > 1.
\end{equation*}
A more interesting problem is to understand if there exists the growing-up solution. We refer to \cite{I-team-Invent} for the case of $\T^2$ and \cite{Chabert} for the case of quantum harmonic oscillator.
\end{remark}

The proof of Theorem \ref{thm1} is based on rescaling argument and I-method developed in \cite{CKSTT2002Almost}. Since the $\T$-direction is not scale-invariant, all the analysis should be done in the rescaled manifold $\R\times\T_{\lambda}$. If $u$ is a solution to $\eqref{NLS-intro}$ posed on $\RT$, then $v = u^{\lambda}(t,x,y)=\lambda^{-1} u(\lambda^{-2}t,\lambda^{-1}x,\lambda^{{-1}}y)$ solves the following equation
\begin{equation}\label{NLS-lambda}
\begin{cases}
iv_t + \Delta_{\rt}v = |v|^2v,\quad (t,x,y) \in \R_t\times\rt \\
v(0,x,y) = v_0(x,y):=\lambda^{-1} u_0(\lambda^{-1}x,\lambda^{{-1}}y),
\end{cases}
\end{equation} 
where $\T_{\lambda} := \lambda \T$.
Let us describe the $I$-method, which consists in  smoothing out the $H^s$ initial data for $0<s<1$, in order to obtain a good local and global solution at $H^1$ level. In our setting, we adapt the I-operator to $\mathcal{F}(Iu)(\xi) := m(\xi) \mathcal{F}(u)(\xi)$, where $\xi \in \R\times\Z_{1/\lambda}$ with $\Z_{1/\lambda} := \frac{1}{\lambda} \Z$ and $m(\xi)$ is a smooth cut off function
\begin{equation}\label{multiplier}
m(\xi) = \begin{cases}
1, & |\xi| \le N,\\
\big(\frac{N}{|\xi|}\big)^{1-s}, & |\xi|\ge 2N.
\end{cases}
\end{equation}
Acting $I$-operator to both side of \eqref{NLS-lambda}, we obtain the modified equation associated to $I$-operator
\begin{align}\label{I-intro}
\begin{cases}
i\partial_t(Iv) + \Delta_{\rt}Iv = I(|v|^2v),\quad (t,x,y) \in \R_t\times\rt \\
Iv(0,x,y) = Iv_0(x,y):=\lambda^{-1} Iu_0(\lambda^{-1}x,\lambda^{{-1}}y).
\end{cases}
\end{align}
Since $Iu^{\lambda}$ is not a solution to $\eqref{NLS-lambda}$, the energy $E(Iu^{\lambda})$ will not be conserved. herefore, it is essential to estimate the increment of energy $E(Iu^{\lambda})$over the time interval during which the local solution exists. This increment estimate is crucial for  establishing the global well-posedness in $H^s$ below the energy space. Recently, Deng-Fan-Yang-Zhao-Zheng\cite{DFYZZ}  successfully removed  the $\varepsilon$-loss of the bilinear estimate on rescaled waveguide $\R\times\T_{\lambda}:=\R\times\lambda\T$. 
Applying this refined bilinear estimate together with the I-method from \cite{CKSTT2002Almost}, they showed the following energy increment 
\begin{equation}\label{Deng}
| E(Iu^{\lambda}(\delta)) - E(Iu^{\lambda}(0)) | \lesssim \lambda^{-\frac{1}{2}}N^{-1+},
\end{equation}
where $\lambda \sim N^{\frac{1-s}{s}}$, $\delta\sim 1$ and $\|Iu^{\lambda}_0\|_{H^1} \lesssim 1$. As a direct consequence, they obtained the global well-posedness for $\eqref{NLS-intro}$ in $H^s$ with $s>\frac{3}{5}$. In order to improve the global well-posedness result to $s > \frac{1}{2}$,  we aim to find a modified energy which enjoys the better decay compared to \eqref{Deng}.

Let us outline   the proof of Theorem \ref{thm1}. The  key point is to construct the good modified energy $E_0(u)$. In \cite{DFYZZ}, it can be observed that the energy increment can be expressed as two terms, where the first term is associated to $\Delta Iu$ and the second term is associated to $I(|u|^2u)$. From their analysis, we can find that the second term will contribute a better decay than the first term, hence the decay estimate is  primarily  determined by the first term. In \cite{CKSTT2008Res}, they first introduce the modified energy by incorporating  a “correction term” into $E(Iu)$. This modification serves to balance both terms effectively. After adding this correction, we observe that the decay of the second term becomes weaker, whereas that of the first term is enhanced. Thus we can expect to show the better decay than \eqref{Deng}.  To be more precisely, we explain and adapt the ideas presented in \cite{CKSTT2008Res} within a semi-periodic framework..
	
Let $k$ be an integer, we set the hyper-plane $\Sigma_k \subset (\R \times \Z_{1/\lambda})^k$ as
\begin{equation*}
\Sigma_k := \{ (\xi_1,\ldots,\xi_k) \in (\R \times \Z_{1/\lambda})^k: \xi_1 + \ldots + \xi_k = 0 \},
\end{equation*}
equipped with the measure which is induced from Lebesgue measure $d\xi_1 \ldots d\xi_{k-1}$. This is obtained by pushing forward under the map
\begin{equation*}
(\xi_1,\ldots,\xi_{k-1}) \mapsto (\xi_1,\ldots,\xi_{k-1},-\xi_1-\ldots-\xi_{k-1}).
\end{equation*} 

We also denote that $\xi_{ijk}=\xi_i+\xi_j+\xi_k$ and $\xi_{ij}=\xi_i+\xi_j$ for $i,j,k\in\{1,2,\cdots,6\}$.

Next, we rewrite the energy $E(Iu)$ to the following form
\begin{align*}
E(Iu)& =\frac{1}{2}\int_{\Sigma_2} |\xi_1||\xi_2|m(\xi_1)m(\xi_2)\widehat{u}(t,\xi_1)\widehat{\overline{u}}(t,\xi_2)(d\xi_1)_\lambda(d\xi_2)_\lambda\notag\\
&\hspace{3ex}+\frac{1}{4(2\pi)^2}\int_{\Sigma_4} m(\xi_1)m(\xi_2)m(\xi_3)m(\xi_4) \widehat{u}(t,\xi_1)\widehat{\overline{u}}(t,\xi_2)\widehat{u}(t,\xi_3)\widehat{\overline{u}}(t,\xi_4).
\end{align*}

To prove Theorem \ref{thm1}, we construct the modified energy $E_0(u)$ by adding the correction term. 
\begin{align}\label{fml-I-energy}
E_0(u)& =\frac{1}{2}\int_{\Sigma_2} |\xi_1||\xi_2|m(\xi_1)m(\xi_2)\widehat{u}(t,\xi_1)\widehat{\overline{u}}(t,\xi_2)(d\xi_1)_\lambda(d\xi_2)_\lambda\notag\\
& \hspace{3ex}+\frac{1}{4(2\pi)^2}\int_{\Sigma_4} \Lambda_4(\xi_1,\xi_2,\xi_3,\xi_4) \widehat{u}(t,\xi_1)\widehat{\overline{u}}(t,\xi_2)\widehat{u}(t,\xi_3)\widehat{\overline{u}}(t,\xi_4),
\end{align}
where
\begin{align}\label{lambda-4}
\Lambda_4(\xi_1,\xi_2,\xi_3,\xi_4)=& \frac{m(\xi_1)^2|\xi_1|^2-m(\xi_2)^2|\xi_2|^2+m(\xi_3)^2|\xi_3|^2-m(\xi_4)^2|\xi_4|^2}{|\xi_1|^2-|\xi_2|^2+|\xi_3|^2-|\xi_4|^2}\notag\\
& \times 1_{S}(\xi_1,\xi_2,\xi_3,\xi_4).
\end{align}
The set $S$ is defined as follows.  Let $N_1,N_2,N_3,N_4$ be dyadic numbers, we set
\begin{align}\label{fml-Res-set}
S(N_1,N_2,N_3,N_4):=&\big\{ (\xi_1,\xi_2,\xi_3,\xi_4) \in \Sigma_4:
|\xi_1| \sim N_1, |\xi_2| \sim N_2, |\xi_3| \sim N_3, |\xi_4| \sim N_4,\nonumber\\
& |\cos \angle(\xi_{12}, \xi_{14})| \geq max\{N_1,N_2,N_3,N_4 \}^{-1} \text{or} \max_{1 \leq j \leq 4} |\xi_j| \ll N \big\},  
\end{align}
where $N$ is as in \eqref{multiplier} and 
\begin{equation}\label{fml-set-S}
S := \bigcup_{(N_1,N_2,N_3,N_4)\in (2^{\mathbb{N}})^4} S(N_1,N_2,N_3,N_4).
\end{equation}

By using the multiplier of $E(Iu)$ and $E_0(u)$, we can deduce that two energies $E_0(u)$ and $E(Iu)$ are comparable, see Proposition \ref{prop-point} below. Therefore, it suffices to prove the energy increment for $E_0(u)$.	

For convenience, we denote $u^\lambda$ by $u$. If $u$ solves \eqref{NLS-lambda}, then $u$ satisfies that
\begin{align*}
u_t=i\Delta_{\RT}u-iu\overline{u}u.
\end{align*}
Taking the Fourier transform to the both side of above equation, we have that
\begin{align}\label{Fourier}
\widehat{u}_t(t,\xi)=-i|\xi|^2\widehat{u}(t,\xi)-\frac{i}{(2\pi)^2}[\widehat{u}*\widehat{\overline{u}}*\widehat{u}](t,\xi),\quad \xi\in\rz.
\end{align}
Combining $\eqref{fml-I-energy}$ and $\eqref{Fourier}$, we obtain
\begin{align*}
& \frac{d}{dt}\big(E_0(u)\big)\\
= & -\frac{i}{2}\int_{\Sigma_2}|\xi_1||\xi_2|(|\xi_1|^2-|\xi_2|^2)m(\xi_1)m(\xi_2)\widehat{u}(t,\xi_1)\widehat{\overline{u}}(t,\xi_2)(d\xi_1)_\lambda (d\xi_2)_\lambda\\
& +\frac{i}{2(2\pi)^2}\int_{\Sigma_4}\left(|\xi_1|^2m(\xi_1)^2-|\xi_2|^2m(\xi_2)^2\right)\widehat{u}(t,\xi_1)\widehat{\overline{u}}(t,\xi_2)\widehat{u}(t,\xi_3)\widehat{\overline{u}}(t,\xi_4)(d\xi_1)_\lambda\cdots(d\xi_4)_\lambda\\
& -\frac{i}{4(2\pi)^2}\int_{\Sigma_4} \Lambda_4(\xi_1,\xi_2,\xi_3,\xi_4) (|\xi_1|^2-|\xi_2|^2+|\xi_3|^2-|\xi_4|^2)\widehat{u}(t,\xi_1)\widehat{\overline{u}}(t,\xi_2)\widehat{u}(t,\xi_3)\widehat{\overline{u}}(t,\xi_4)(d\xi_1)_\lambda\cdots(d\xi_4)_\lambda\\
& -\frac{i}{4(2\pi)^4}\int_{\Sigma_6}\Lambda_6(\xi_1,\xi_2,\xi_3,\xi_4,\xi_5,\xi_6)\widehat{u}(t,\xi_1)\widehat{\overline{u}}(t,\xi_2)\widehat{u}(t,\xi_3)\widehat{\overline{u}}(t,\xi_4)\widehat{u}(t,\xi_5)\widehat{\overline{u}}(t,\xi_6)(d\xi_1)_\lambda\cdots(d\xi_6)_\lambda,
\end{align*}
where
\begin{align*}
\Lambda_6(\xi_1,\xi_2,\xi_3,\xi_4,\xi_5,\xi_6) = & \Lambda_4(\xi_{123},\xi_4,\xi_5,\xi_6) - \Lambda_4(\xi_1,\xi_{234},\xi_5,\xi_6) + \Lambda_4(\xi_{1},\xi_2,\xi_{345},\xi_6) - \Lambda_4(\xi_{1},\xi_2,\xi_3,\xi_{456}).
\end{align*}
Notice that the first two terms is identical to zero since $\xi_1+\xi_2=0$. Thus, we obtain
\begin{align}\label{E0}
	&\hspace{3ex} \frac{d}{dt}E_0(u)(t)\notag \\
	&=  -\frac{i}{4(2\pi)^2}\int_{\Sigma_4} A(\xi_1,\xi_2,\xi_3,\xi_4)\widehat{u}(t,\xi_1)\widehat{\overline{u}}(t,\xi_2)\widehat{u}(t,\xi_3)\widehat{\overline{u}}(t,\xi_4)(d\xi_1)_\lambda\cdots(d\xi_4)_\lambda\notag\\
	& \hspace{3ex}-\frac{i}{4(2\pi)^4}\int_{\Sigma_6}\Lambda_6(\xi_1,\xi_2,\xi_3,\xi_4,\xi_5,\xi_6)\widehat{u}(t,\xi_1)\widehat{\overline{u}}(t,\xi_2)\widehat{u}(t,\xi_3)\widehat{\overline{u}}(t,\xi_4)\widehat{u}(t,\xi_5)\widehat{\overline{u}}(t,\xi_6)(d\xi_1)_\lambda\cdots(d\xi_6)_\lambda,
\end{align}
where
\begin{align*}
  A(\xi_1 , \xi_2 , \xi_3 , \xi_4) := & \big( m(\xi_1)|\xi_1|^2 - m(\xi_2)|\xi_2|^2 + m(\xi_3)|\xi_3|^2 - m(\xi_4)|\xi_4|^2 \big)  1_{ S^c}(\xi_1 , \xi_2 , \xi_3 , \xi_4).
\end{align*}

	

By the definition of $E(Iu)$ and $E_0(u)$,  $\eqref{fml-I-energy}$ and $\eqref{E0}$, Theorem \ref{thm1}  is reduced to  the following two propositions.

\begin{proposition}[Pointwise bound]\label{prop-point}
If $\|u_0\|_{L^2(\R\times\T_{\lambda})} \leq 1$ and $E(Iu_0) \leq 1$, and $Iu\in C([0,\delta],H^1(\R\times\T_{\lambda}))$ is a solution to \eqref{I-intro} on a time interval $[0,\delta]$ where the solution is locally well-posed. Then we have
\begin{equation}\label{pointbound}
\left|E_0(u)-E(Iu)\right|\lesssim \frac{1}{N^{1-}}\|Iu\|_{H^{1}(\rt)}^4.
\end{equation} 
\end{proposition}
\begin{proposition}[Energy-increment for $E_0(u)$]\label{Prop1}
 Let $\delta\sim1$ and $Iu\in C([0,\delta],H^1(\rt))$ be the solution to \eqref{I-intro} satisfy $\|Iu_0\|_{H^1(\rt)}\lesssim1$, then we have
\begin{align}\label{prop-es}
|E_0(u(\delta))-E_0(u_0)|\lesssim N^{-2+}\big(\|Iu\|_{X^{1,b}(\rt)}^4+\|Iu\|_{X^{1,b}(\rt)}^6\big),
\end{align}
where $X^{s,b}$ denotes the Bourgian space, which will be defined in Section $\ref{Pre}$.
\end{proposition}
By the definition of $E_0(u)$, the proof of Proposition \ref{Prop1} can be reduced to prove the following two propositions. 
\begin{proposition}[Sextilinear estimate]\label{sextilinear-estimate}
Let $\delta,u_0,u$ and $Iu$ be the same as Proposition \ref{Prop1}. Then we have
\begin{align}\label{eq:e-6}
& \left|\int_{0}^\delta\int_{\Sigma_6}\Lambda_6(\xi_1,\xi_2,\xi_3,\xi_4,\xi_5,\xi_6)\widehat{u}(t,\xi_1)\widehat{\overline{u}}(t,\xi_2)\widehat{u}(t,\xi_3)\widehat{\overline{u}}(t,\xi_4)\widehat{u}(t,\xi_5)\widehat{\overline{u}}(t,\xi_6)\,(d\xi_1)_\lambda\cdots\,(d\xi_6)_\lambda\,dt\right|\notag\\
\lesssim  & \frac{1}{N^{2-}}\|Iu\|_{X^{1,1/2+}}^6.
\end{align} 
\end{proposition}
\begin{proposition}[Quadrilinear-estimate]\label{quadrilinear-prop}
Let $\delta,u_0,u$ and $Iu$ be the same as Proposition \ref{Prop1}. Then we have
\begin{equation}\label{eq-6}
\left|\int_{0}^\delta\int_{\Sigma_4} A(\xi_1,\xi_2,\xi_3,\xi_4)\widehat{u}(t,\xi_1)\widehat{\overline{u}}(t,\xi_2)\widehat{u}(t,\xi_3)\widehat{\overline{u}}(t,\xi_4)\,(d\xi_1)_\lambda\cdots\,(d\xi_4)_\lambda\,dt\right|\lesssim \frac{1}{N^{2-}}\|Iu\|_{X^{1,1/2+}}^4.
\end{equation}
\end{proposition}
For the pointwise bound \eqref{pointbound} that describe the difference between $E(Iu)$ and $E_0(u)$, we can reduce it to the following pointwise bound, by using the definition of $E_0(u)$ and $E(Iu)$
\begin{equation*}
\bigg|\int_{\Sigma_4}\big(\Lambda_4(\xi_1,\xi_2,\xi_3,\xi_4) - m(\xi_1)m(\xi_2)m(\xi_3)m(\xi_4)\big)\widehat{u}(\xi_1)\widehat{\bar{u}}(\xi_1) \widehat{u}(\xi_3)\widehat{\bar{u}}(\xi_4)\bigg|\lesssim N^{-1+}\|Iu\|_{X^{1,b}(\rt)}^4.
\end{equation*}
We will use the Newton-Leibniz formula to establish the multiplier estimate. We refer to Section \ref{Local} for more details. 
	
	For the sextilinear estimate \eqref{eq:e-6}, we will apply pointwise bound for $|\Lambda_4(\xi_1,\xi_2,\xi_3,\xi_4)|$ and $L^4$ Strichartz estimate to treat it in Section \ref{Six}.
	
	The quadrilinear estimate is more complicated. To handle this term, we establish the following refined bilinear Strichartz estimate with angular truncation:
\begin{theorem}[Angularly refined bilinear estimate in $\R\times\T_\lambda$]\label{bil-intro}
		Let $0 < \theta \ll 1 < M$ and $1 <N_2 \le N_1$. Suppose $\phi_{N_1}$ and $\phi_{N_2}$ are two functions which are frequency localized around $N_1$ and $N_2$ respectively. We define a bilinear function as
		\begin{equation*}
			F(t,z) := \int_{(\rz)^2} e^{-it (|\xi_1|^2 + |\xi_2|^2) -iz \cdot (\xi_1 + \xi_2) } 1_{|\mu_1 - \mu_2|\gtrsim M} 1_{| cos\angle (\xi_1 , \xi_2) |\le \theta}\widehat{\phi_{N_1}}(\xi_1) \widehat{\phi_{N_2}}(\xi_2) (d\xi_1)_{\lambda} (d\xi_2)_{\lambda},
		\end{equation*}
		where $z=(x,y)\in\R\times\T_{\lambda}$ and $\xi_j=(\mu_j,\eta_j)\in\rz$ with $\Z_{1/\lambda}=\tfrac1\lambda\Z$ for $j=1,2$.
		Then we have
		\begin{equation}\label{fml-int-es-bi}
			\| F \|_{L^2([0,1] \times \R \times \mathbb{T}_{\lambda})} \lesssim \bigg( \frac{1}{\lambda M} + \theta \bigg)^{\frac{1}{2}} \|\phi_{N_1}\|_{L^2} \|\phi_{N_2}\|_{L^2}.
		\end{equation}
	\end{theorem}
\begin{remark}
 Such type refined bilinear estimate was initially established in \cite{CKSTT2008Res}. Recently, Takaoka \cite{Takaoka} showed the above bilinear estimate \eqref{fml-int-es-bi} with $\lambda=1$.  In this proposition, we extend the result of \cite{Takaoka} to the case of $\lambda\geq1$. We also notice that when $\lambda$ approaches infinity, our estimate goes back to that in \cite{CKSTT2008Res} without any loss of derivatives. 
\end{remark}

	This paper is organized as follows. In Section \ref{Pre}, we collect the basic properties of functional spaces and Littlewood-Paley theory on waveguide setting. We also give the definition of I-operator and construction of the modified energy $E_0(u)$ associated to $Iu$. In Section \ref{Local}, we prove the local well-posedness of $\eqref{I-intro}$ and control the differences of $E(Iu)$ and $E_0(u)$. The proof of main theorem can be reduced to the increment estimate of $E_0(u)$. In Section \ref{Six} and Section \ref{Four}, we prove the sextilinear   and quadrilinear estimate which imply the energy increment.
	
	\subsection{Notations}
	We denote $A\lesssim B$, if there exist constant $C > 0$ such that $A \le C B$. We say that $A \sim B$ if both $A \lesssim B$ and $B \lesssim A$ hold. We write $A \ll B$ if for some small constant $0 < c\ll1$ such that $A \le c B$. We use the notation  $\langle a \rangle := (1 + |a|^2)^{\frac{1}{2}}$ and $ a\pm := a \pm \varepsilon $.
	\section{Preliminaries}\label{Pre}
	In this section, we will briefly recall the functional space and the associated harmonic analysis tools on the waveguide manifold. We also give the definition and properties  of  $I$-operator and the modified energy associated to $Iu$.
	
	\subsection{Functional spaces and Littlewood-Paley theory on waveguide manifold}
	
	Throughout this paper, we denote the Lebesgue space $L_t^pL_z^q(I\times\R\times\T)$ with the norm 
	\begin{align*}
		\|u\|_{L_t^pL_z^q(I\times\R\times\T)}=\bigg(\int_{I}\Big(\int_{\R\times\T}|u(t,z)|^qdz\Big)^\frac{p}{q}\bigg)^\frac1p,
	\end{align*}
	where $I$ is the time interval and $z=(x,y)\in\R\times\T$.
	
	We define the Fourier transform on $\R\times\T$ as follows:
	\begin{align*}
		(\mathcal{F}f)(\xi)=\int_{\R\times\T}f(z)e^{-iz\cdot\xi}dz,
	\end{align*}
	where $\xi=(\mu,\eta)\in\R\times\Z$. Meanwhile, we have the Fourier inverse formula
	\begin{align}
		f(z)=\sum_{\eta\in\Z}\int_{\R}(\mathcal{F}f)(\xi)e^{iz\cdot\xi}d\mu.
	\end{align}
	For convenience, we denote that
	\begin{align}\label{notation-1}
		\int_{\R\times\Z}f(\xi)d\xi\stackrel{\triangle}{=}\sum_{\eta\in\Z}\int_{\R}f(
		\mu,\eta)d\mu.
	\end{align}
	In the rest of the paper, we shall use the  notation \eqref{notation-1} since the discrete sum can be understood as the integral with discrete measure.
	
	Moreover, we write the linear Schr\"odinger operator as the following 
	\begin{align*}
		\mathcal{F}(e^{it\Delta_{\R\times\T}}f)(\xi)=e^{-it|\xi|^2}(\mathcal{F}f)(\xi).
	\end{align*} 
	Now, we introduce the Littlewood-Paley projections. First, we denote a smooth cut-off function $\eta_1(\xi)$ by the following
	\begin{align*}
		\eta_1(\xi)=\begin{cases}
			1,&|\xi|\le\frac34,\\
			0,&|\xi|\ge\frac83.
		\end{cases}
	\end{align*}
	Let $N\in2^{\Bbb{Z}}$ be a dyadic number, then we define the Littlewood-Paley operator  $P_{\le N}$ and $P_N$ by
	\begin{align*}
		\mathcal{F}(P_{\le N}f)(\xi)=\eta_1(\frac{\xi}{N})\mathcal{F}f(\xi), \quad \xi\in\R\times\Z,
	\end{align*}
	and 
	\begin{align*}
		P_{N}f=P_{\le N}f-P_{\le \frac{N}{2}}f.
	\end{align*}
	For arbitrary $a>0$, we define
	\begin{align*}
		P_{\le a}:=\sum_{N\le a}P_Nf,\quad P_{>a}f=\sum_{a<N}P_Nf.
	\end{align*}
	
	Scaling symmetry plays a crucial role  in the global well-posedness of nonlinear Schr\"odinger equation on $\R^d$. However, when  performing  the rescaling on a manifold, it is essential to carefully consider the spatial domain, since the metric and curvature of the rescaled manifold may be altered. Now we give the definition of the Fourier transformation and Littlewood-Paley projection on the rescaled waveguide $\R\times\T_\lambda$ as well.
	
	We  denote the Fourier transform on $\R \times \T_{\lambda}$ by the following
	\begin{equation*}
		\hat{f}(\xi)= (\mathcal{F} f)(\xi)= \int_{\mathbb{R} \times \mathbb{T}_\lambda}f(z)e^{-2\pi i z\cdot \xi}dz,
	\end{equation*} 
	where $z = (x,y) \in \R \times \T_{\lambda}$ and $\xi = (\mu , \eta) \in \R \times \Z_{1/\lambda}$. We also use the following notation for convenience
	\begin{equation}\label{notation-2}
		\int_{\R\times\Z_{1/\lambda}} f(\xi)(d \xi)_{\lambda}:=\frac{1}{\lambda} \sum_{\eta \in\mathbb{Z}_{1/\lambda}} \int_{ \mathbb{R}}f(\mu,\eta) d\mu.
	\end{equation}
	By the direct calculus, we have the following  lemma. 
	\begin{lemma}\label{lem-Fourier-proper}
		The following estimates hold for functions defined on $\R\times\T_{\lambda}$:
		\begin{enumerate}
			\item Fourier inverse formula
			\begin{equation}
				f(x,y)= \frac{1}{\lambda} \sum_{\eta \in \mathbb{Z}_{1/\lambda}} \int_{ \mathbb{R}} \widehat{f}(\mu,\eta)e^{2\pi i(x\mu + y\eta)}d\mu.
			\end{equation}
			\item Plancherel identity
			\begin{equation}
				\|f\|_{L^2\left(\mathbb{R} \times \mathbb{T}_\lambda\right)}=\|\hat{f}\|_{L^2\left(\mathbb{R} \times \mathbb{Z}_{1/\lambda}\right)  }.  
			\end{equation}
			\item Convolution and multiplication under Fourier transform
			\begin{equation}
				\widehat{f g}(\xi)=\hat{f} * \hat{g}(\xi):=\int_{\mathbb{R} \times \mathbb{Z}_{1/\lambda}} \hat{f}\left(\xi-\xi^\prime\right) \hat{g}\left(\xi^\prime\right)\left(d \xi^\prime\right)_\lambda.
			\end{equation}
		\end{enumerate}
	\end{lemma}
	Using the inverse Fourier formula, the solution $U_\lambda(t)u_0$ to the linear Schr\"odinger equation 
	\begin{align}\label{fml-linear}
		i\partial_tu(t,z)+\Delta_{\R\times\T_{\lambda}}u(t,z)=0,\quad u|_{t=0}=u_0,\quad z\in\R\times\T_\lambda
	\end{align}
	can be rewritten as
	\begin{align*}
		U_{\lambda}u_0 (t,z) = \int_{\R\times\Z_{1/\lambda}} e^{2\pi i z \cdot \xi - (2\pi |\xi|)^2 it} \widehat{u_0} (\xi) (d\xi)_{\lambda}.
	\end{align*}
	
	Following from \cite{Bour1993,Bour1998Refine} and let $s > 0$ and $b \in \R$. We define the $X^{s,b}$ norm associated to $\Delta_{\R\times\T_{\lambda}}$ as
	\begin{equation}
		\| u \|_{X^{s,b}(\R\times\R\times\T_{\lambda})} := \| \langle \tau - |\xi|^2 \rangle^b \langle \xi \rangle^s ( \mathcal{F}_{t,z} u ) \|_{L^2_{\tau , \xi}(\R\times\R\times\Z_{1/\lambda})},
	\end{equation}
	where $u(t,z)$ is a function defined on $\R\times\R\times\T_{\lambda}$. We denote the space-time Fourier transform of $u$ by $\mathcal{F}_{t,z} u(\tau,\xi)$, where $(\tau,\xi) \in \R \times \R \times \Z_{1/\lambda}$. 
	Therefore, we can define the space-time conormal space $X^{s,b}$, which contains all functions with finite $X^{s,b}$ norms.
	\begin{definition}[Bourgain sapce]
		Let $s\geqslant0$ and $b\in\Bbb{R}$,
		\begin{align*}
			X^{s,b}(\Bbb{R}\times \rt)=\left\{u\in\mathcal{ S}(\Bbb{R},L^2(\rt)):\Vert u\Vert_{X^{s,b}(\Bbb{R}\times \rt)}<\infty\right\}.
		\end{align*}
		For $s\leqslant0$, we can define the negative Bourgain space by duality. Moreover, we define $X^{\infty,b}$ as $ X^{\infty,b}=\bigcap\limits_{s\in\Bbb{R}}X^{s,b}$. 
		
		For $T>0$, we can define the restricted Bourgain space equipped with the norm
		\begin{align*}
			\Vert u\Vert_{X_T^{s,b}}=\inf\left\{\Vert w\Vert_{X^{s,b}(\Bbb{R}\times \rt)}<\infty,\hspace{1ex}w|_{[-T,T]\times \rt}=u\right\}.
		\end{align*}
	\end{definition}
	Then we introduce Strichartz estimate obtained by Barron\cite{Bar2021GloStr} on rescaled waveguide
	\begin{lemma}[Strichartz estimate, \cite{Bar2021GloStr}]
		Let $4 \le p \le \infty$, $\alpha(p) = 1 - \frac{4}{p}$,  then we have
		\begin{equation}\label{fml-L4}
			\|u\|_{L^p(\R\times\R\times\T_{\lambda})} \lesssim  \|u\|_{X^{\alpha(p),\frac{1}{2}+}(\R\times\T_{\lambda})}.
		\end{equation}
	\end{lemma}
	The key role of Bourgain space is that we can inherit bilinear estimate from linear solution. To be more specifically,  we have the following transfer principle as in \cite{Bour1998Refine,Burq2005bilinear}.
	\begin{lemma}
		Let $b > \frac{1}{2}$ and $N_1,N_2$ be any dyadic numbers, then the following two statements are equivalent: 
		\begin{enumerate}
			\item For any space function $f_1,f_2 \in L^2(\R\times\T_{\lambda})$, we can find constant  $C(\lambda , N_1 , N_2)>0$ such that
			\begin{align}
				\| U_{\lambda}(t) P_{N_1}f_1 U_{\lambda}(t) P_{N_2}f_2 &\|_{L^2_{t,z} ([0,1]\times \R \times \T_{\lambda}) }\nonumber \\
				&\lesssim C(\lambda , N_1 , N_2) \|f_1\|_{L^2(\R\times\T_{\lambda})}\|f_2\|_{L^2(\R\times\T_{\lambda})}.
			\end{align} 
			\item For any space-time function $u_1,u_2 \in X^{0,b}(\R \times \R \times \T_{\lambda})$, we can find constant $C(\lambda , N_1 , N_2)>0$ such that
			\begin{align}
				\| P_{N_1}u_1  P_{N_2}u_2 &\|_{L^2_{t,z} ([0,1]\times \R \times \T_{\lambda}) }\nonumber \\ 
				&\lesssim C(\lambda , N_1 , N_2) \|u_1\|_{X^{0,b}(\R\times\R\times\T_{\lambda})}\|u_2\|_{X^{0,b}(\R\times\R\times\T_{\lambda})}.
			\end{align}
		\end{enumerate}	
	\end{lemma}

\section{The setting of I-method and  proof of Theorem \ref{thm1}}\label{Local}
As mentioned in introduction, the I-method played the important role in studying low regularity global well-posedness of the nonlinear dispersive equations. Since the lack of conservation law at the $H^s$ level with $0<s<1$, the key ingredient of I-operator is to smooth the $H^s$ data to $H^1$. In our setting, it can be defined as follows:
\begin{definition}
	Let $N$ be dyadic number, we set $\mathcal{F}(Iu)(\xi) := m(\xi) \mathcal{F}(u)(\xi)$ where $\xi \in \rz$, $m(\xi)$ is a smooth function such that
	\begin{equation}
		m(\xi) = \begin{cases}
			1, & |\xi| \le N,\\
			\big(\frac{N}{|\xi|}\big)^{1-s}, & |\xi|\ge 2N.
		\end{cases}
	\end{equation}
\end{definition}

In the following, we shall prove Theorem \ref{thm1} by applying the I-method and under Proposition \ref{Prop1}. First, we introduce the scaling of two equation posed in different geometries. If $u(t,z)$ solves Cauchy problem $\eqref{NLS-intro}$ on $\R\times\T$, then the function $u^{\lambda}(t,z)$
\begin{align*}
	u^{\lambda}(t,z)=\frac{1}{\lambda}u\big(\frac{t}{\lambda^2},\frac{x}{\lambda},\frac{y}{\lambda})
\end{align*}
solves the nonlinear Schr\"{o}dinger equation posed on $\R\times\T_{\lambda}$,
\begin{align}\label{fml-NLS-lambda}
	\begin{cases}
		i\partial_tu+\Delta_{\R\times\T_{\lambda}} u= \left|u\right|^2u, &(t,z)\in\mathbb{R}\times \R\times\T_{\lambda},\\
		u(0,z)=u_0(z)\in H^s(\R\times\T_{\lambda}).
	\end{cases}
\end{align}
In the rest of this paper, we replace $u^{\lambda}$ by $u$ for convenience.
\subsection{The local well-posedness of the I-system}
We impose the $I$ to the both sides of the Cauchy problem $\eqref{fml-NLS-lambda}$, and denote $v$ as $Iu$, which satisfies
\begin{align}\label{I}
	\begin{cases}
		i\partial_tv+\Delta_{\R\times\T_{\lambda}} v=I(\left|u\right|^2u)=I(\left|I^{-1}v\right|^2I^{-1}v),&(t,z)\in\Bbb{R}\times \R\times\T_{\lambda},\\
		v(0)=v_0(z)=Iu_0.
	\end{cases}
\end{align}
The usual energy conservation law fails in the case of $Iu$, but we can prove the energy of $Iu$ is not increasing very fast in the small time interval $[0,\delta]$ and then we call it the almost conservation law. We define the pseudo-energy as follows
\begin{align*}
	E(Iu(t))&=\int_{\R\times\T_{\lambda}}\Big(\frac{1}{2}\left|\nabla
	Iu\right|^2+\frac{1}{4}\left|Iu\right|^4\Big)dz\\
	&=\frac{1}{2}\Vert Iu\Vert_{\dot{H}^1}^2+\frac{1}{4}\Vert Iu\Vert_{L^4}^4.
\end{align*}
For convenience, we also regard it as ``energy'' of $Iu$.

In order to prove the local well-posedness of $Iu$, we need to control the term $I(\left|u\right|^2u)$, which is given by the following lemma.
\begin{lemma}[Nonlinear estimate, \cite{Takaoka}]\label{I-estimate}
	Let $(b,b')\in\Bbb{R}^2$ and $\delta > 0$ satisfy
	\begin{align*}
		0<b'<\frac{1}{2}<b,\quad b+b'<1.
	\end{align*}
	For  $u\in X_\delta^{s,b}$ with $s>\frac{1}{2}$, there exists a constant $C>0$ such that
	\begin{align*}
		\Vert I(\left|u\right|^2u)\Vert_{X_\delta^{1,b'}}\leqslant C\Vert Iu\Vert^3_{X_\delta^{1,b}}.
	\end{align*}
\end{lemma}

Indeed, Lemma \ref{I-estimate} implies the local well-posedness for the Cauchy problem \eqref{I}.
\begin{proposition}[Local well-posedness]\label{Local-M}
	Suppose $\Vert Iu_0\Vert_{H^1(\R\times\T_{\lambda})}\lesssim1$. Then, there exists $\delta\sim1$ such that  the Cauchy problem $\eqref{I}$ is local well-posed on the  small time interval $[0,\delta]$. Furthermore, let $b>\frac{1}{2}$, we have
	\begin{equation*}
		\| Iu \|_{X^{1,b}_\delta} \lesssim \| Iu_0 \|_{H^1(\R\times\T_{\lambda})} \lesssim 1.
	\end{equation*}
\end{proposition}

\subsection{Point-wise estimate}
In this subsection, we aim to show the pointwise bound, i.e. Proposition \ref{prop-point}   by using the definition of $E_0(u)$ and $E(Iu)$.
This pointwise bound give a good control on difference between  two types of energies: $|E(Iu)-E_0(u)|$. Therefore, it remains to obtain the almost conservation law for the modified energy $E_0(u)$. 

By the definition of $E_0(u)$, we only need to estimate
\begin{equation}
	\bigg|\Lambda_4(\xi_1,\xi_2,\xi_3,\xi_4) - m(\xi_1)m(\xi_2)m(\xi_3)m(\xi_4) \bigg|\lesssim N_1\min\{m(N_1),m(N_2),m(N_3),m(N_4)\}^2.
\end{equation}
By analyzing the symbol, we will prove the pointwise estimate.
\begin{proposition}[Point-wise bound]\label{prop-pointwise}  
	Suppose that $Iu \in H^1(\R\times\T_{\lambda})$ is the solution to equation $\eqref{I}$, then for every $t > 0$, we have
	\begin{equation}\label{fixed-bound2}
		| E(Iu(t)) - E_0(u(t)) | \lesssim  N^{-1+} \| Iu \|_{H^1_z(\R\times\T_{\lambda})}^4.
	\end{equation}
\end{proposition}
To prove the above proposition, we need the following estimates.
	\begin{lemma}[\cite{CKSTT2008Res,Takaoka}]\label{lemma-pointwise-multi}  For any $(\xi_1,\xi_2,\xi_3,\xi_4) \in \Sigma_4$, We have
		\begin{equation*}
			\big|m(\xi_1)^2|\xi_1|^2 - m(\xi_2)^2|\xi_2|^2 + m(\xi_3)^2|\xi_3|^2 - m(\xi_4)^2|\xi_4|^2\big| \lesssim \min( m(\xi_1),m(\xi_2),m(\xi_3),m(\xi_4))^2 |\xi_{12}| |\xi_{14}|,
		\end{equation*}
		and 
		\begin{equation*}
			|\Lambda_4(\xi_1,\xi_2,\xi_3,\xi_4)| \lesssim N_1\min(m(N_1),m(N_2),m(N_3),m(N_4))^2.
		\end{equation*}
		where $\xi_{12} = \xi_1 + \xi_2$ and $\xi_{14} = \xi_1 + \xi_4$.
	\end{lemma}
	We also need the following pointwise estimate for $\Lambda_4$.
	\begin{corollary}\label{cor-bound}  For any $(\xi_1, \xi_2, \xi_3,
		\xi_4) \in \Sigma_4$ and $|\xi_j| \sim N_j, j=1,2,3,4$ where $N_j\in 2^{\mathbb{N}}$. Then, we have
		\begin{equation*}
			|\Lambda_4(\xi_1,\xi_2,\xi_3,\xi_4)| \lesssim N_1\min(m(N_1),m(N_2),m(N_3),m(N_4))^2.
		\end{equation*}
	\end{corollary}
	Using Lemma \ref{lemma-pointwise-multi} and Corollary \ref{cor-bound}, we can  prove Proposition \ref{prop-point} as in \cite{CKSTT2008Res}.

\subsection{The proof of Theorem \ref{thm1}}
Let $u,u_0$ be the same as in Theorem \ref{thm1}. For the fixed $T>0$, we define $u:[0,T]\times\R\times\T\to \Bbb{C}$ and $A:=1+\|u(0)\|_{H^s_{z}(\RT)}$. Suppose that  $\lambda$ is  the scaling parameter to be chosen later, we  define $u^{\lambda}$ as 
	\begin{equation*}
		u^{\lambda}(t,z) := \frac{1}{\lambda} u(\frac{t}{\lambda^2}, \frac{x}{\lambda},\frac{y}{\lambda}),\quad z=(x,y)\in\R\times\T_{\lambda}.
	\end{equation*}
	Therefore, we  find that $u^{\lambda}:[0,\lambda^2T]\times\R\times\T_{\lambda}\to\Bbb{C}$ is the solution to $\eqref{fml-NLS-lambda}$. Then by the Gagliardo-Nirenberg inequality \cite{Hebey}, we have
	\begin{equation*}
		\|Iu^{\lambda}\|_{L^4(\R\times\T_{\lambda})}^4\lesssim \|Iu^{\lambda}\|^2_{\dot{H}^{1}(\R\times\T_{\lambda})}\|Iu^{\lambda}\|^2_{L^2(\R\times\T_{\lambda})}.
	\end{equation*}
	Hence, the  energy quantity of $Iu^\lambda$ is totally controlled by its kinetic energy. Indeed, we have		
	\begin{equation*}
		E(Iu^{\lambda}(0)) \lesssim \|Iu^{\lambda}(0)\|^2_{\dot H^{1}(\rt)} \lesssim N^{2-2s}\lambda^{-2s} A^4 .
	\end{equation*}
	
	Taking $\lambda = c(s,A) N^{\frac{1-s}{s}}$ with small constant $c(s,A)$, we have the uniform bound
	\begin{align}\label{fml-gwp-EIu0}
		E(Iu^\lambda(0))\le\frac{1}{3}.
	\end{align} 
	By the mass conservation law, we also have
	\begin{equation}\label{fml-gwp-Iu-L2}
		\|Iu^{\lambda}(t)\|_{L^2(\rt)} \le A.
	\end{equation}
	From $\eqref{fml-I-energy}$, $\eqref{fml-gwp-Iu-L2}$ and the pointwise bound Proposition \ref{prop-point}, we get
	\begin{equation}\label{fml-gwp-pointwise}
		| E(Iu^\lambda(t)) - E_0(u^{\lambda}(t)) | \lesssim N^{-1+} (A^4 + E(Iu^{\lambda}(t))^2).
	\end{equation}
	Taking $\lambda=c(s,A)N^{\frac{1-s}s}$ where $c(s,A)$ is chosen sufficiently small and $N$ is chosen sufficiently large,	we claim that  the following estimate holds for any $0 < t \le \lambda^2 T$ 
	\begin{equation}\label{fml-gwp-claim}
		E(Iu^{\lambda}(t)) < \frac{2}{3}.
	\end{equation}
	With this in hand, we obtain
	\begin{equation*}
		\|u^{\lambda}(\lambda^2 t)\|_{H^s(\rt)} \lesssim \|Iu^{\lambda}(\lambda^2 t)\|_{H^1(\RT)} \lesssim A + E(Iu^{\lambda}(\lambda^2 t)) \lesssim 1 + A,\quad\forall 0\leqslant t\leqslant T.
	\end{equation*}
	We finally obtain
	\begin{equation*}
		\|u(t)\|_{H^s(\RT)} \lesssim \lambda^s (1 + A),\quad 0\leqslant t\leqslant T,
	\end{equation*}
	which implies that $\eqref{NLS-intro}$ is well-posed in space $C([0,T];H^s(\RT))$. Hence, we have proved Theorem \ref{thm1}.
	\begin{remark}
		Indeed, one can obtain the polynomial growth of Sobolev norms
		\begin{align*}
			\|u(T)\|_{H^s(\RT)}\leqslant (1+\|u(0)\|_{H^s})^{C}T^{\frac{s(1-s)}{4s-2}+},
		\end{align*}
		where $C$  depends on $s$. 
	\end{remark}

\section{The proof of sextilinear estimate}\label{Six}
In this subsection, as an application of $L^4$ Strichartz estimate, we will prove Proposition \ref{sextilinear-estimate}.
	
	\begin{proof}
		The left-hand side of \eqref{eq:e-6} can be expressed as
		\begin{equation*}
			\left|\int_0^{\delta} \int_{\Sigma_6} \Lambda_6(\xi_1,\xi_2,\xi_3,\xi_4,\xi_5,\xi_6)
			\hat u(t,\xi_1) \ldots \hat{\overline{u}}(t,\xi_6) \ dt\right|.
		\end{equation*}
		On the region $max\{|\xi_1|, \cdots, |\xi_6| \} \le N/3$, we claim that $\Lambda_6(\xi_1,\xi_2,\xi_3,\xi_4,\xi_5,\xi_6)  =0$. Indeed, by definition we find that
		\begin{align*}
			\Lambda_6(\xi_1,\xi_2,\xi_3,\xi_4,\xi_5,\xi_6) = & \Lambda_4(\xi_{123},\xi_4,\xi_5,\xi_6) - \Lambda_4(\xi_1,\xi_{234},\xi_5,\xi_6) \\
			& + \Lambda_4(\xi_{1},\xi_2,\xi_{345},\xi_6) - \Lambda_4(\xi_{1},\xi_2,\xi_3,\xi_{456}) = 0,
		\end{align*} 
		via the fact $|\xi_{123}|,|\xi_{234}|,|\xi_{345}|,|\xi_{456}| \le N$. Therefore, it remains to consider the region $\max\{|\xi_1|,\ldots,|\xi_6|\} \geq N/3$. Proposition \eqref{sextilinear-estimate} is reduced to the following estimate
		\begin{equation}\label{fml-six-whole}
			I:= \left|\int_0^{\delta} \int_{\Sigma_6 } 1_{\{\max(|\xi_1|,\ldots,|\xi_6|) \geq N/3\}}\Lambda_6(\xi_1,\xi_2,\xi_3,\xi_4,\xi_5,\xi_6)
			\hat u(t,\xi_1) \ldots \hat{\overline{u}}(t,\xi_6) \ dt\right| \lesssim N^{-2+}\|Iu\|_{X^{1,b}}^6.
		\end{equation}
		Using Littlewood-Paley decomposition to $\hat{u}(t,\xi_j)$, we can deduce $\eqref{fml-six-whole}$ to the following estimate
		\begin{align}\label{fml-sexty-target}
			I(N_1,N_2,N_3,N_4,N_5,N_6)&:= \bigg|\int_0^{\delta} \int_{\Sigma_6 } 1_{\{\max(|\xi_1|,\ldots,|\xi_6|) \geq N/3\}}\\
			&\hspace{12ex} \times \Lambda_6(\xi_1,\xi_2,\xi_3,\xi_4,\xi_5,\xi_6)
			\hat u(t,\xi_1) \ldots \hat{\overline{u}}(t,\xi_6) \ dt\bigg|\nonumber\\
			&\lesssim (N_1N_2N_3N_4N_5N_6)^{0-}N^{-2+}\|Iu\|_{X^{1,b}}^6\nonumber,
		\end{align}
		where $\hat u(t,\xi_j)$ are compactly supported in the interval $[N_j,2N_j]$ with $N_j\in 2^{\mathbb{N}}$ for $j = 1,2,3,4,5,6$. To treat the right-hand side of \eqref{fml-sexty-target}, we will follow the strategy in \cite{CKSTT2008Res}.
		
		Without loss of generality, assuming that $N_1\ge N_2 \ge N_3\ge N_4 \ge N_5\ge N_6\geq1$ and we will prove $\eqref{fml-sexty-target}$ for different cases. If at least three frequencies are no less than $N$, that is $N_1\ge N_2 \ge N_3 \gtrsim N$.
		Because the $X^{s,b}$ norm uses the spacetime Fourier transform, we will be forced for technical reasons to write the left-hand side 
		in terms of the spacetime Fourier transform.  Indeed, this left-hand side becomes
		\begin{align*} 
			\left|\int_{\R^6}  \int_{\Sigma_6} 
			\hat 1_{[0,t_0]}(\tau_0) 1_{\max(|\xi_1|,\ldots,|\xi_6|) \geq N/3} \Lambda_6(\xi_1,\xi_2,\xi_3,\xi_4,\xi_5,\xi_6)\tilde u(\tau_1,\xi_1)
			\ldots \tilde {\overline{u}}(\tau_6,\xi_6)(d\xi_1)_\lambda\cdots(d\xi_6)_\lambda d\tau_1 \ldots d\tau_6\right|
		\end{align*}
		where $\tilde{u}$ is the space-times Fourier transform of $u$ and $\tau_0 := -\tau_1 - \ldots -\tau_6$.  Using  the following bound
		\begin{equation}\label{t0-fourier}
			\left|\hat 1_{[0,t_0]}(\tau_0)\right| \lesssim \langle \tau_0 \rangle^{-1},
		\end{equation}
		we can control this quantity by
		\[ \lesssim \sup_{\xi_ j,j=1,\cdots,6}|\Lambda_6(\xi_1,\xi_2,\cdots,\xi_6)|
		\int_{\R^6} \int_{\Sigma_6} 
		\langle \tau_0 \rangle^{-1} 1_{\max(|\xi_1|,\ldots,|\xi_6|) \geq N/3} 
		|\tilde u(\tau_1,\xi_1)|
		\ldots |\tilde{\overline{u}}(\tau_6,\xi_6)|(d\xi_1)_{\lambda}\cdots(d\xi_6)_{\lambda} d\tau_1 \ldots d\tau_6.
		\]
		
		Without loss of generality, we will not distinguish $u$ and its conjugation a when using the $X^{s,b}$ norm.
		so it suffices to show the estimate
		\begin{align*} 
		\sup_{\xi_j,j=1,\cdots,6}	\big|\Lambda_6(\xi_1,\xi_2,\xi_3,\xi_4,\xi_5,\xi_6)\big|\cdot	\int_{\R^6} & \int_{\Sigma_6} 
			\langle \tau_0 \rangle^{-1} 1_{\max(|\xi_1|,\cdots,|\xi_6|) \geq N/3} 
			\prod_{j=1}^6 |\tilde u_j(\tau_j,\xi_j)|(d\xi_1)_{\lambda}\cdots(d\xi_6)_{\lambda} 
			d\tau_1 \ldots d\tau_6\\
			& \lesssim N^{-2+}N_1^{0-} \prod_{j=1}^6 \|Iu_j\|_{ X^{1,b}}.\end{align*}
		
		Inspired by \cite{CKSTT2008Res},  the $X^{s,b}$ norm depends only on the size (i.e., magnitude) of the spacetime Fourier transform. This allows us to take all $\tilde{u}_j$ as non-negative, and consequently, the absolute value signs can be omitted.
		The left-hand side can now be written using spacetime convolutions as
		\begin{equation}\label{mn}
			\int_\R \langle \tau \rangle^{-1}\tilde u_1 * \cdots *  \tilde u_6(\tau,0)\ d\tau.
		\end{equation}
		Unfortunately, the function $\langle \tau \rangle^{-1}$ is not integrable. To overcome this difficulty,  we introduce the logarithmic weight $w(\tau) := 1 + \log^2 \langle \tau \rangle$ inspired by \cite{CKSTT2008Res}, then $\langle \tau \rangle^{-1} w^{-1}$ is integrable.  Also, from the elementary estimate $w(\tau_1 + \ldots + \tau_6) \lesssim w(\tau_1) \ldots w(\tau_6)$ we have the pointwise bound
		\[ \tilde  u_1 * \cdots * \tilde  u_6 \lesssim w^{-1} [ (w \tilde u_1) * \ldots (w \tilde u_6) ].
		\]
		Thus, we have the following estimate
		$$ \eqref{mn} \lesssim \| (w\tilde  u_1) * \ldots * (w\tilde  u_6) \|_{L^\infty_{\tau,\xi}}.$$
		Thus it will suffice to show that
		\[\sup_{\xi_j,j=1,\cdots,6}\big|\Lambda_6(\xi_1,\xi_2,\xi_3,\xi_4,\xi_5,\xi_6)\big|\cdot\| (w \tilde u_1) * \ldots * (w \tilde u_6) \|_{L^\infty_{\tau,\xi}}
		\lesssim N^{-2+} N_1^{0-} \prod_{j=1}^6 m(N_j) N_j \| u_j \|_{ X^{0,\frac{1}{2}+}}.
		\]
		If $v_j$ denotes the function with spacetime Fourier transform $ v_j = w  u_j$, one easily verifies that
		$$ \| v_j \|_{ X^{0,b-}} \lesssim \log^2 (1 + N_1) \| u_j \|_{ X^{0,b}}.$$
		We shall take $b=\frac{1}{2}+\varepsilon$. By  Hausd\"orff-Young's inequality, it remains to prove
		\begin{align}\label{remain} \big|\Lambda_6(\xi_1,\xi_2,\xi_3,\xi_4,\xi_5,\xi_6)\big|\cdot\|  v_1 \cdots v_6 \|_{L^1_{t,z}}
			\lesssim N^{-2+} N_1^{0-} \prod_{j=1}^6 m(N_j) N_j \| v_j \|_{ X^{0,b-}}.
		\end{align}

		We apply Corollary \ref{cor-bound}, which gives upper bound to $|\Lambda_6(\xi_1,\xi_2,\xi_3,\xi_4,\xi_5,\xi_6)|$. Meanwhile, taking $L^4$ Strichartz estimate $\eqref{fml-L4}$ to functions $v_j$ with $j = 1,2,3,4$, Sobolev embedding $X^{1+,\frac{1}{2}+} \hookrightarrow L_{t,z}^{\infty}$ to functions $v_5,v_6$ and apply the H\"older inequality to estimate $\|v_1\cdots v_6\|_{L_{\tau,\xi}^1}$. 
		Therefore, we obtain
		\begin{align*}
			&\hspace{6ex}\mbox{(LHS) of }\eqref{remain}\\
			&\lesssim N_1m(N_4)^2 N_5^{1+}N_6^{1+} \prod_{j=1}^6  \|u_{N_j}\|_{X^{0,\frac{1}{2}+}_\delta(\rt)} \\
			&\lesssim N_1m(N_4)^2 N_5^{0+}N_6^{0+} \prod_{j=1}^6 m(N_j)^{-1} N_j^{-1}  \|Iu_{N_j}\|_{X^{1,b}_\delta(\rt)} \\
			&\lesssim \frac{m(N_4)^2}{m(N_5)m(N_6)} \frac{1}{m(N_1)N_1^{1/2}m(N_2)N_2^{1/2}} \frac{1}{m(N_3)N_3^{1-}m(N_4)N_4^{1-}}  \prod_{j=1}^6\|Iu_{N_j}\|_{X^{1,\frac{1}{2}+}_\delta(\rt)},
		\end{align*}
		where we use the fact that $N_1\sim N_2$.
		It is clear that if $\alpha\geq\frac{1}{2}>1-s$, we have for $j=1,2,3$
		\begin{equation*}
			m(N_j)N_j^{\alpha} = m(N_j)N_j^{\frac{1}{2}} N_j^{\alpha - \frac{1}{2}} \gtrsim N^{\alpha}.
		\end{equation*} 
		Notice that $m$ is a non-increasing function, thus we have that
		\begin{align*}
			\frac{m(N_4)^2}{m(N_5)m(N_6)}\lesssim1.
		\end{align*}
		Thus, we obtain
		\begin{equation*}
			\frac{m(N_4)^2}{m(N_5)m(N_6)} \frac{1}{m(N_1)N_1^{1/2}m(N_2)N_2^{1/2}} \frac{1}{m(N_3)N_3m(N_4)N_4}\lesssim N^{-2} \lesssim N^{-2+} N_1^{0-}.
		\end{equation*}
		
		Now, we turn to the case that only two frequencies are larger than $N$, that is $N_1 \ge N_2 \gtrsim N \gg N_3$. In this case,   we need a improved upper bound for $|\Lambda_6(\xi_1,\xi_2,\xi_3,\xi_4,\xi_5,\xi_6)|$. By the definition, we have
		\begin{align}\label{fml-six-Lambda6}
			|\Lambda_6(\xi_1,\xi_2,\xi_3,\xi_4,\xi_5,\xi_6)| \lesssim & |\Lambda_4(\xi_{123},\xi_4,\xi_5,\xi_{6})| + |\Lambda_4(\xi_{1},\xi_{234},\xi_5,\xi_{6})|\nonumber\\
			& + |\Lambda_4(\xi_{1},\xi_2,\xi_{345},\xi_{6})| + |\Lambda_4(\xi_{1},\xi_2,\xi_3,\xi_{456})|.
		\end{align}
		It is sufficient to consider the pointwise bound for $|\Lambda_4(\xi_1^{\prime} , \xi_2^{\prime} , \xi_3^{\prime} , \xi_4^{\prime})|$ with $\xi_1^{\prime} + \xi_2^{\prime} + \xi_3^{\prime} + \xi_4^{\prime} = 0$. Thanks to $\eqref{fml-six-Lambda6}$, we only need to consider the following cases
		\begin{enumerate}
			\item \textbf{Case 1:} $|\xi_j^{\prime}| \ll N(\xi_1^{\prime} = \xi_{123},\xi_2^{\prime} = \xi_4,\xi_3^{\prime} = \xi_5,\xi_4^{\prime} = \xi_6)$, \\
			\item \textbf{Case 2:} $|\xi_1^{\prime}|,|\xi_2^{\prime}| \gtrsim N$, $|\xi_3^{\prime}|,|\xi_4^{\prime}| \ll N(\text{Related to other three terms in right-hand side of }\eqref{fml-six-Lambda6})$.
		\end{enumerate}

		\noindent\textbf{Case 1:}$|\xi_j^{\prime}| \ll N$. This case is trivial, since $\Lambda_4(\xi_1^{\prime} , \xi_2^{\prime} , \xi_3^{\prime} , \xi_4^{\prime}) = 0$. 
		\\
		\\
		\noindent\textbf{Case 2:}$|\xi_1^{\prime}|,|\xi_2^{\prime}| \gtrsim N,|\xi_3^{\prime}|,|\xi_4^{\prime}| \ll N$. In this case, we always have $|\xi_1^{\prime}| \sim |\xi_2^{\prime}| \gtrsim N \gg |\xi_3^{\prime}| + |\xi_4^{\prime}| $. By Newton-Leibniz and standard argument in \cite{CKSTT2008Res}, we can find
		\begin{align}\label{fml-six-claim}
			|\Lambda_{4}(\xi_1^{\prime}&,\xi_2^{\prime},\xi_3^{\prime},\xi_4^{\prime})|\nonumber\\
			\lesssim & max \{m(\xi_1^{\prime}),m(\xi_2^{\prime}),m(\xi_3^{\prime}),m(\xi_4^{\prime})\}^2 \frac{|\xi_3^{\prime}| + |\xi_4^{\prime}|}{N_1^{-1}(|\xi_1^{\prime}| + |\xi_2^{\prime}|)}.
		\end{align}
		Hence, $|\Lambda_{4}(\xi_1^{\prime},\xi_2^{\prime},\xi_3^{\prime},\xi_4^{\prime})|$ enjoys the bound
		\begin{align}\label{fml-six-Lam4-bound}
			|\Lambda_{4}(\xi_1^{\prime}&,\xi_2^{\prime},\xi_3^{\prime},\xi_4^{\prime})|
			\lesssim max \{m(\xi_1^{\prime}),m(\xi_2^{\prime}),m(\xi_3^{\prime}),m(\xi_4^{\prime})\}^2 N_3.
		\end{align}
		Arguing like the previous case, we have the following estimate by combining  $\eqref{fml-six-Lam4-bound}$, H\"older's inequality, $L^4$ Strichartz estimate with Sobolev embedding $X^{1+,\frac{1}{2}+} \hookrightarrow L^{\infty}$
		\begin{align*}
			I(N_1,N_2,N_3,N_4,N_5,N_6)
			\lesssim &m(N_5)m(N_6) N_3 N_5^{1+}N_6^{1+} \prod_{j=1}^6  \|u_{N_j}\|_{X^{0,\frac{1}{2}+}_\delta(\rt)} \\
			\lesssim& m(N_5)m(N_6) N_3 N_5^{1+}N_6^{1+} \prod_{j=1}^6 m(N_j)^{-1} N_j^{-1}  \|Iu_{N_j}\|_{X^{1,\frac{1}{2}+}_\delta(\rt)} \\
			\lesssim& \frac{1}{m(N_1)N_1m(N_2)N_2^{1-}} \frac{1}{m(N_3)N_3^{0+}m(N_4)N_4^{1-}}  \prod_{j=1}^6\|Iu_{N_j}\|_{X^{1,\frac{1}{2}+}_\delta(\rt)}\\
			\lesssim& N^{-2+}N_1^{0-} \prod_{j=1}^6\|Iu_{N_j}\|_{X^{1,b}_\delta(\rt)} \lesssim N^{-2+}\prod_{j=1}^{6}\|Iu_j\|_{X^{1,\frac{1}{2}+}_\delta(\rt)}.
		\end{align*}
		Thus we have proved $\eqref{fml-sexty-target}$. 
	\end{proof}
\section{The proof of quadrilinear estimate}\label{Four}
In this section, we will prove the quadrilinear estimate by establishing the angularly refined bilinear Strichartz estimates. The bilinear Strichartz estimate was first developed in \cite{Bour1998Refine}
\begin{equation}
	\| e^{it\Delta} u_0 e^{it\Delta} v_0 \|_{L^2_{t,x}(\R\times\R^2)} \lesssim \big(\frac{N_2}{N_1}\big)^{\frac{1}{2}} \|u_0\|_{L^2(\R^2)} \|v_0\|_{L^2(\R^2)} ,
\end{equation}
where $u_0,v_0$ are frequency-localized functions around the dyadic number $N_1\ge N_2$. Later, Colliander-Keel-Staffilani-Takaoka-Tao \cite{CKSTT2008Res} established the following angularly refined bilinear Strichartz estimate. 
\begin{lemma}[Angularly refined bilinear estimate,\cite{CKSTT2008Res}]\label{Tao}Let $0<N_2\leq N_1$ and $0<\theta\ll\frac{N_2}{N_1}$. Then for any $f_1,f_2\in L^2(\R^2)$ with spatial frequencies $N_1$ and $N_2$ respectively, then the bilinear function
	\begin{align*}
		F(t,x)=\int_{\R^2}\int_{\R^2}e^{-it(|\xi_1|^2+|\xi_2|^2)}e^{ix\cdot(\xi_1+\xi_2)}1_{|\cos\angle(\xi_1,\xi_2)|\leq\theta}\widehat{f_1}(\xi_1)\widehat{f_2}(\xi_2)d\xi_1d\xi_2
	\end{align*}
	enjoys the bound
	\begin{align*}
		\|F\|_{L_{t,x}^2(\R\times\R^2)}\lesssim\theta^\frac12\|f_1\|_{L^2(\R^2)}\|f_2\|_{L^2(\R^2)}.
	\end{align*}
\end{lemma}
Recently, Takaoka \cite{Takaoka} established the similar angularly refined bilinear estimate on waveguide $\R\times\T$. 
\begin{lemma}[Angularly refined bilinear estimate on $\R\times\T$,\cite{Takaoka}]\label{Takaoka}
	Let $0<\theta\ll1<M$ and $1<N_2\leq N_1$. Suppose that $f_{N_1},f_{N_2}\in L^2(\R\times\T)$ are functions with spatial frequencies $N_1$ and $N_2$ respectively. Then we define the space-time function as
	\begin{align*}
		F(t,z)=\int e^{-it(|\xi_1|^2+|\xi_2|^2)}e^{iz\cdot(\xi_1+\xi_2)}1_{|\mu_1-\mu_2|\gtrsim M}1_{|\cos\angle(\xi_1,\xi_2)|\leq \theta}\widehat{f_{N_1}}(\xi_1)\widehat{f_{N_2}}(\xi_2)d\xi_1d\xi_2
	\end{align*}
	where $z=(x,y)\in\R\times\T$ and $\xi_j=(\mu_j,\eta_j)\in\R\times\Z$ for $j=1,2$. Then it satisfies
	\begin{align*}
		\|F\|_{L_{t,z}^2([0,1]\times\R\times\T)}\lesssim\frac{\langle\theta N_1\rangle^{1/2}}{M^{1/2}}\|f_{N_1}\|_{L^2(\RT)}\|f_{N_2}\|_{L^2(\RT)}.
	\end{align*}
\end{lemma}

\subsection{Angularly refined bilinear estimate on the rescaled waveguide $\R\times\T_{\lambda}$}
 We will prove Theorem \ref{bil-intro} by counting arguments in Fourier space, taking advantage of the product measure. 
	
	\begin{proof}[Proof of Theorem \ref{bil-intro} ]
	   Indeed, by performing the following angular decomposition
		\begin{equation*}
			\phi_{N_1} = \sum_{l_1} \phi_{N_1,l_1}, \quad \phi_{N_2} = \sum_{l_2} \phi_{N_2,l_2}, 
		\end{equation*}
		where $l_1,l_2$ are integers between $0$ and $2\pi /\theta$, and the Fourier transform of $\phi_{N_j,l_j}$ is supported in the sector $\{ \xi_j \in \rz | |arg(\xi_j) - l_j\theta| \lesssim \theta \}$. Hence, we can rewrite the bilinear form $F(t,z)$ as
		\begin{align}\label{reduce-2-angular}
			F(t,z) = \sum_{l_1,l_2} \int_{(\rz)^2} & e^{-it (|\xi_1|^2 + |\xi_2|^2) + iz \cdot (\xi_1 + \xi_2) } 1_{|\mu_1 - \mu_2|\gtrsim M}\nonumber \\ 
			&\times 1_{| cos\angle (\xi_1 , \xi_2) |\le \theta}\widehat{\phi}_{N_1,l_1}(\xi_1) \widehat{\phi}_{N_2,l_2}(\xi_2) (d\xi_1)_{\lambda} (d\xi_2)_{\lambda}.
		\end{align}
        
        Without loss of generality, we can assume the support of $\widehat{\phi_{N_j}},j = 1,2$ is contained in two angular sectors 
		\begin{equation*}
			\{ \xi_j | |arg(\xi_j) - l_j \theta| \lesssim \theta \} \subset \rz,
		\end{equation*}
		where $l_j$ are integers satisfying $1 \le l_j \le \frac{2\pi}{\theta}$. Thus, we need to consider the following bilinear form
		\begin{align}
        \label{angular-reduce}
			F(t,z) := & \int_{(\rz)^2} e^{-it (|\xi_1|^2 + |\xi_2|^2) + iz \cdot (\xi_1 + \xi_2) } 1_{|\mu_1 - \mu_2|\gtrsim M} \widehat{\phi_{N_1}}(\xi_1) \widehat{\phi_{N_2}}(\xi_2) \notag\\
			& \times 1_{\arg(\xi_1) = l_1 \theta + O(\theta)} 1_{\arg(\xi_2) = l_2 \theta + O(\theta)} (d\xi_1)_{\lambda} (d\xi_2)_{\lambda}.
		\end{align}

		Following the standard arguments in \cite{Silva2007Global,DFYZZ}, we have
		\begin{align}\label{fml-qua-bi-A}
			\|F\|_{L_{t,z}^2([0,1]\times\rt )}\lesssim \sup\limits_{(\tau , \xi)\in (\R \times \rz)} |K(\tau , \xi)|  \|\phi_{N_1}\|_{L^2(\rt)} \|\phi_{N_2}\|_{L^2(\rt)},
		\end{align}		
		where recalling that
		\begin{align}\label{fml-M-unmbercont}
			K(\tau , \xi) = \Big\{ & \xi_2 = (\mu_2 , \eta_2) \in \R \times \Z_{1/\lambda} |\nonumber \\ & arg(\xi_2) = l_1\theta + O(\theta), arg(\xi - \xi_2) = l_2 \theta + O(\theta), |\xi_2| \sim N_2, |\xi - \xi_2| \sim N_1, \nonumber \\ & | \mu/2 - \mu_2| \gtrsim M,  |\tau - |\xi_2|^2 - |\xi - \xi_2|^2 | \le 1  \Big\},
		\end{align}
		with $\xi = (\mu,\eta) \in \rz$ and $\tau \in \R$. 
		Hence the estimate of the bilinear form $F(t,z)$ is reduced  to a measure estimate problem, i.e.
		\begin{align}\label{fml-claim-coun}
			\sup\limits_{(\tau , \xi)\in (\R \times \rz)} |K(\tau,\xi)| \lesssim \frac{1}{\lambda M}+\theta.
		\end{align}
		Here, the factor $\frac{1}{\lambda}$ comes from the normalized measure $\R\times\T_{\lambda}$.
		We first claim that $K(\tau,\xi)$ is totally contained in an annulus $S(r)$. We claim that $S(r)$ can be expressed as: 
		\begin{equation}\label{fml-four-Sr}
			S(r) := \left\{ \xi_2 : |\xi_2 - \frac{\xi}{2}| = r + O(1/r) \right\},
		\end{equation}
		with 
		\begin{equation*}
			r = \frac{\sqrt{2\tau-\xi^2}}{2}.
		\end{equation*}
		Moreover, $r$ enjoys the bound $M\lesssim r \lesssim N_1$. By direct computation, we get
		\begin{equation*}
			\bigg|\tau - \frac{\xi^2}{2} - 2(\frac{\xi}{2} - \xi_2)^2 \bigg| = |\tau - |\xi_2|^2 - |\xi - \xi_2|^2 | \le 1,
		\end{equation*}
		which implies
		\begin{equation*}
			\frac{2\tau - \xi^2 - 4}{4} \le (\frac{\xi}{2} - \xi_2)^2 \le \frac{2\tau - \xi^2 + 4}{4}.
		\end{equation*}
		Thus, using Taylor expansion, we find that the width of $S(r)$ is $O(1/r)$.
		
		From $\tau \sim |\xi_2|^2 + |\xi - \xi_2|^2$, we can find
		\begin{equation*}
			N_1^2 \gtrsim 2\tau - \xi^2 \sim 4(\frac{\xi}{2} - \xi_2)^2 \gtrsim M^2,
		\end{equation*}
		and the claim follows.
		
		We are devoted to prove \eqref{fml-claim-coun} when $N_2 \ll N_1$ and $N_1 \sim N_2$ respectively.
		
		For the \textbf{case $N_2 \ll N_1$}, we prove $\eqref{fml-claim-coun}$ in the following three sub-cases.
		
		\textbf{Sub-case one: }$|\mu| \sim N_1$. If $\xi_2 , \tilde{\xi}_2 \in K(\tau,\xi)$ satisfy $\eta_2 = \tilde{\eta}_2$, we have
		\begin{equation*}
			| \tau - |\xi_2|^2 - |\xi - \xi_2|^2 | \le 1,\text{and}\quad | \tau - |\tilde{\xi}_2|^2 - |\xi - \tilde{\xi}_2|^2 | \le 1, 
		\end{equation*}
		which implies
		\begin{equation*}
			2|(\xi_2 - \tilde{\xi}_2)\cdot(\xi_2 + \tilde{\xi}_2 - \xi)| = | |\xi_2|^2 - |\tilde{\xi}_2|^2 + |\xi - \xi_2|^2 - |\xi - \tilde{\xi}_2|^2 | \lesssim 1.
		\end{equation*}
		Notice that $\xi_2 - \tilde{\xi}_2 = (\mu_2 - \tilde{\mu}_2,0)$, hence  we get
		\begin{equation*}
			|(\mu_2 - \tilde{\mu}_2)\cdot(\mu_2 + \tilde{\mu}_2 - \mu)| \lesssim 1,
		\end{equation*}
		notice that $|\mu_2 + \tilde{\mu}_2 - \mu| \gtrsim N_1$,
		Thus, for each fixed $k\in \Z$, we have
		\begin{equation*}
			\bigg| K(\tau,\xi) \bigcap \R \times \{k/\lambda\} \bigg| \lesssim \frac{1}{N_1}.
		\end{equation*}
		Now, we devote to estimate the number of $k$ such that $K(\tau,\xi) \bigcap \R \times \{k/\lambda\}$ is non-empty. 
		
		To obtain this, we consider rectangles of $O(1/\lambda) \times O(1/r)$. Next, we cover the annulus $S(r)$ with the above rectangles and split the angular $arg(\xi_2-\xi) = l_2\theta + O(\theta)$ into $(1 + \lambda \theta r)$ pieces, which can be expressed as
		\begin{equation*}
			arg(\xi_2-\xi) \in l_2\theta + \big[ \frac{k}{\lambda r}, \frac{k+1}{\lambda r} \big], |k|\lesssim (1 + \lambda \theta r).
		\end{equation*} 
		\begin{center}
			\begin{tikzpicture}[scale=1]
			\draw [thick, domain = 0: 90] plot ({2.2*cos(\x)-3} , {2.2*sin(\x)})  ;
			\draw [thick, domain = 90: 180] plot ({2.2*cos(\x)-3} , {2.2*sin(\x)})  ;
			\draw [thick, domain = 180: 270] plot ({2.2*cos(\x)-3} , {2.2*sin(\x)})  ;
			\draw [thick, domain = 270: 360] plot ({2.2*cos(\x)-3} , {2.2*sin(\x)})  ;
			\draw [thick, domain = 0: 90] plot ({2*cos(\x)-3} , {2*sin(\x)})  ;
			\draw [thick, domain = 90: 180] plot ({2*cos(\x)-3} , {2*sin(\x)})  ;
			\draw [thick, domain = 180: 270] plot ({2*cos(\x)-3} , {2*sin(\x)})  ;
			\draw [thick, domain = 270: 360] plot ({2*cos(\x)-3} , {2*sin(\x)})  ;
			
			\draw [thick] (-3,0) -- (-3,-2.2) ;
			\fill (-3,0) circle (1pt) ; 
			\fill (-3,-2) circle (1pt) ; 
			\fill (-3,-2.2) circle (1pt) ; 
			\draw (-3,0) node [left] {O} ;
			\draw (-3,-2) node [left=0.2cm,above] {A} ;
			\draw (-3,-2.2) node [left,below] {B} ;
			
			\draw (-3,-3) node {$|OA|=r,|AB|=O(1/r)$} ;
			
			\draw [-] (-1.1,-0.3) rectangle (-0.7,0.3) ;
			\draw [rotate around={14:(-3,0)}] (-1.1,-0.3) rectangle (-0.7,0.3) ;
			\draw [rotate around={-14:(-3,0)}] (-1.1,-0.3) rectangle (-0.7,0.3) ;	
			\draw [thick] (-3,0) -- (-1.1,-0.3) ;
			\draw [thick] (-3,0) -- (-1.1,0.3) ;
			\draw [rotate around={-14:(-3,0)}][thick] (-3,0) -- (-1.1,-0.3) ;
			\draw [rotate around={14:(-3,0)}][thick] (-3,0) -- (-1.1,0.3) ;	
				
			\fill [rotate around={35:(-3,0)}] (-2.2,0) circle (0.7pt) ;
			\fill [rotate around={45:(-3,0)}] (-2.2,0) circle (0.7pt) ;
			\fill [rotate around={55:(-3,0)}] (-2.2,0) circle (0.7pt) ;
	
			\fill [rotate around={-35:(-3,0)}] (-2.2,0) circle (0.7pt) ;
			\fill [rotate around={-45:(-3,0)}] (-2.2,0) circle (0.7pt) ;
			\fill [rotate around={-55:(-3,0)}] (-2.2,0) circle (0.7pt) ;
			
			\draw [line width=0.5pt, double distance=2pt,
			arrows = {-Latex[length=0pt 3 0]}] (-0.3,0) -- (2,0) ;
			
			\draw [thick] (2.5,-1.5) rectangle (3.5,1.5) ;
			\draw [<->] (2.5,-1.6) -- (3.5,-1.6) ;
			\draw (3,-1.6) node[below] {$1/r$} ;
			\draw [<->] (3.6,-1.5) -- (3.6,1.5) ;
			\draw (3.6,0) node[right] {$\frac{1}{\lambda}$} ;
			
			\path (0,-4) node[xshift = 0.5cm](caption) {Figure 1: Decompose annulus $S(r)$ into rectangles.} ;
			\end{tikzpicture}
		\end{center}
		
		Each angular piece $arg(\xi_2-\xi) \in l_2\theta + [k/(\lambda r), (k+1)/(\lambda r)]$ intersects at most three rectangles we have mentioned above. Hence, we can find $K(\tau,\xi)$ is contained in the part of annulus $S(r)$, and the arc length of this part has the bound $O(1 + \lambda \theta r)$. Consequently, there number of $k$'s such that $K(\tau,\xi) \bigcap \R \times \{k/\lambda\}$ is non-empty can be bounded by $1 + \lambda N_1 \theta$. Therefore, we summary the first case as follows
		\begin{equation*}
			|K(\tau,\xi)| \le 3|K(\tau,\xi)\cap \R\times\{k/\lambda\}| \le  \frac{3(1+\lambda\theta N_1)}{\lambda N_1} \lesssim\frac{1}{\lambda M}+\theta.
		\end{equation*}
		
		\textbf{Sub-case Two: $|\eta| \sim N_1, \lambda \gtrsim N_1$}, if $\xi_2=(\mu_2,\eta_2) , \tilde{\xi}_2=(\tilde{\mu_2},\tilde{\eta_2}) \in K(\tau,\xi)$. By exchanging the role of $\mu$ and $\eta$, we fix $\mu\in\R$ first. Then it holds
        \begin{align*}
            \#\big\{k\in\Z:(\mu,\frac{k}{\lambda}\in K(\tau,\xi)\big\}\leq\frac{\lambda}{N_1}.
        \end{align*}
        . Arguing like the previous case, we have $K(\tau,\xi)$ is contained in the part of $S(r)$ with arc length less than $N_1\theta$. Therefore, the length of $\mu$ is at most $N_1\theta$. By the direct calculus, we can obtain that 
		\begin{equation*}
			|K(\tau,\xi)|\lesssim \frac{1}{\lambda}\int_{\R} \sum_{\eta_2} |K(\tau,\mu,\eta)| d\mu_2 \lesssim \frac{1}{\lambda}\int_{\R} \frac{\lambda}{N_1}\mathbf{1}_{\{|\mu_2|\lesssim N_1\theta\}}(\mu_2)d\mu_2 \lesssim \theta.
		\end{equation*}
		
		Finally, it yields that
		\begin{equation*}
			|K(\tau,\xi)| \lesssim\frac{1}{\lambda M}+\theta.
		\end{equation*}
		
		\textbf{Sub-case Three: $|\eta| \sim N_1, \lambda \ll N_1$}. Without loss of generality, we assume that $\eta > 0$ and $\eta \sim N_1$. For fix $\eta_2 = \frac{k}{\lambda} \in \Z_{1/\lambda}$, we define 
		\begin{equation*}
			C_k(\tau , \xi) := \bigg\{ \xi_2 \in \R \times \Z_{1/\lambda} : \bigg| \big| \mu_2 - \frac{\mu}{2}  \big|^2 - \tau_k \bigg| \le \frac{1}{2}, \eta_2 = \frac{k}{\lambda} \bigg\},
		\end{equation*}
		where 
		\begin{equation*}
			\tau_k = \frac{\tau}{2} - \big| \frac{k}{\lambda} - \frac{\eta}{2} \big|^2 - \frac{|\xi|^2}{4}.
		\end{equation*}
		Then we can find 
		\begin{align*}
			K(\tau,\xi) \subset \bigg( \bigcup_{k \in \Z} C_k \bigg) \bigcap \{ \xi_2 | |\xi_2| \sim N_2  \} \bigcap S(r) .
		\end{align*}
		Since $|\frac{k}{\lambda}|\leq |\xi_2| \sim N_2$ and $| | \mu_2 - \frac{\mu}{2}  |^2 - \tau_k | \le \frac{1}{2}$, we only need to consider all $k$'s that satisfy $|k| \lesssim \lambda N_2$ and $\tau_k>-1/2$ to estimate $|K(\tau,\xi)|$. The direct calculus yields
		\begin{equation*}
			\tau_{k + 1} - \tau_k \sim \frac{N_1}{\lambda} \gg 1,
		\end{equation*}
		which shows that $\{\tau_k\}$ is is an almost arithmetic sequence with a common difference $N_1/\lambda \gg 1$. Thus, there exist at most finite $k$'s such that $|\tau_k| \le 1$, and we can find
		\begin{equation*}
			\bigg| \bigg( \bigcup_{|\tau_k|\le 1} C_k \bigg) \bigcap \{ \xi_2 | |\xi_2| \sim N_2  \} \bigcap S(r) \bigg| \lesssim \sup\limits_{|\tau_k|\le 1} \big| C_k \bigcap S(r) \big|.
		\end{equation*}
		We assert that for any $k\in \Z$
		\begin{equation}\label{fml-quad-line-measure}
			\bigg| C_k \bigcap K(\tau,\xi) \bigg| \lesssim \frac{1}{\lambda M}.
		\end{equation}
		For $k \in \Z$ such that $\big\{ \xi_2=(\mu_2,\eta_2) \in \rz: \eta_2 = \frac{k}{\lambda} \big\} \bigcap K(\tau,\xi)$ is not empty, we get
		\begin{equation*}
			\bigg| \frac{k}{\lambda} - \frac{\eta}{2} \bigg|^2 = \bigg| \xi_2 - \frac{\xi}{2} \bigg|^2 - \bigg| \mu_2 - \frac{\mu}{2} \bigg|^2 \lesssim r^2 - M^2.
		\end{equation*}
		Therefore, the distance from the line $\R\times \{ k/\lambda \}$ to the center of annulus $(\frac{\mu}{2},\frac{\eta}{2})$ is less than $(r^2 - M^2)^{\frac{1}{2}}$. Hence, the worst-case occurs when 
		\begin{equation*}
			\bigg| \frac{k}{\lambda} - \frac{\eta}{2} \bigg| \sim (r^2 - M^2)^{\frac{1}{2}}.
		\end{equation*}
		For a fixed $k$, we denote the point where $\R\times \{ k/\lambda \}$ intersects the $\Z_{1/\lambda}$-axis as $E$, and the points where $\R\times \{ k/\lambda \}$ intersects annulus $S(r)$ as $A$ and $C$. Connecting $OA$ and let it intersect $S(r)$ at point $B$. From $B$, draw a tangent to the boundary of $S(r)$ (perpendicular to $OB$) intersecting $\R\times \{ k/\lambda \}$ at $D$. The process can be illustrated as follows 
		\begin{center}
			\begin{tikzpicture}[scale=1] 
				\draw[->] (-1,0) -- (3,0) node[anchor=north] {$x$};
				\draw[->] (0,-1) -- (0,3)  node[anchor=east] {$y$};
				\draw (0,0) node[below left]{O};
				
				\draw [thick, domain = 35: 75] plot ({2*cos(\x)} , {2*sin(\x)})  ;
				\draw [thick, domain = 35: 75] plot ({1.5*cos(\x)} , {1.5*sin(\x)})  ;
				
				\draw [-] (0,0) -- (1,2)  ;
				
				\draw [-] (-0.4,1.342) -- (2,1.342) ;
				
				\draw [-] (0.894,1.788) -- (1.82,1.342) ;
				

				\draw (0.670,1.342) coordinate (A) ;
				\draw (0.894,1.788) coordinate (B) ;
				\draw (1.82,1.342) coordinate (D) ;
				\draw pic [draw,black,thick,angle radius=1.3mm] {right angle = A--B--D};
				
				\draw (0,1.342) node[left] {E} circle (0.03)[fill=gray!90] ;
				\draw (0.670,1.342) node[left] {A} circle (0.03)[fill=gray!90];
				\draw (0.894,1.788) node[above] {B} circle (0.03)[fill=gray!90];
				\draw (1.82,1.342) node[above] {D} circle (0.03)[fill=gray!90];
				\draw (1.483,1.342) node[above] {C} circle (0.03)[fill=gray!90];
				
				\draw (0,0.1) node[below right,xshift=0cm,yshift=0cm] {};
				\draw (0,-0.1) node[below,xshift=0cm,yshift=0cm] {} ;
				\draw (0.670,1.442) node[above,right,xshift=0cm,yshift=-0.1cm] {};
				
				\draw (-1,1) node[left] {$|OA|=r,|AB|=O(\frac{1}{r})$} ;
				
				\path (0,-1.3) node[xshift = 0.5cm](caption) {Figure 2: $\R\times \{ k/\lambda \}$ cross through the annulus $S(r)$.} ;
			\end{tikzpicture}
		\end{center}
		Here $O$ denote the center of annulus $S(r)$, $x$ and $y$ directions are parallel with $\R$-axis and $\Z_{1/\lambda}$-axis respectively. From Figure 1, we can see that $\angle ABD$ is a right angle which implies the triangle $AOE$ is similar to that of $ABD$. Therefore, we get 
		\begin{equation*}
			\big| \R \times \{ k/\lambda \} \cap K(\tau,\xi) \big| =\frac{1}{\lambda}|AC| < \frac{1}{\lambda}|AD| = \frac{1}{\lambda}\frac{|AO|\cdot|AB|}{|AE|} \lesssim \frac{1}{\lambda M},
		\end{equation*}
		then $\eqref{fml-quad-line-measure}$ follows. Therefore, we get
		\begin{equation*}
			\bigg| \bigg( \bigcup_{k:|\tau_k|\le 1} C_k \bigg) \bigcap \{ \xi_2 | |\xi_2| \sim N_2  \} \bigcap S(r) \bigg| \lesssim \frac{1}{\lambda M}.
		\end{equation*}
		
		It remains to consider when $\tau_k > 1$, for these $k$'s,
		$$
		C_k=\left\{\xi_2 \in \mathbb{R} \times \mathbb{Z}_{1 / \lambda}:\left(\tau_k-\frac{1}{2}\right)^{\frac{1}{2}} \leq\left|\mu_2-\frac{\mu}{2}\right| \leq\left(\tau_k+\frac{1}{2}\right)^{\frac{1}{2}},\eta_2=\frac{k}{\lambda}\right\},
		$$
		thus we can directly estimate
		$$
		\left|C_k\right| \sim \lambda^{-1} \tau_k^{-\frac{1}{2}}.
		$$

		Let
		\begin{equation*}
			S_{ang}(\theta) := \{ \xi_2 | |arg(\xi_2) - l_1\theta| \lesssim \theta, |arg(\xi - \xi_2) - l_2 \theta| \lesssim \theta \},
		\end{equation*}
		then we define
		\begin{equation*}
			k_0 := \min \{ k \in \Z : C_k \bigcap \{ \xi_2 | |\xi_2| \sim N_2  \} \bigcap S(r) \bigcap S_{ang}(\theta) \neq\emptyset,\tau_k>1,|k|\lesssim \lambda N_2\},
		\end{equation*}
		\begin{equation*}
			k_1 := \max \{ k \in \Z : C_k \bigcap \{ \xi_2 | |\xi_2| \sim N_2  \} \bigcap S(r) \bigcap S_{ang}(\theta)\neq\emptyset,\tau_k>1,|k|\lesssim \lambda N_2 \}.
		\end{equation*}
		Using the similar strategy as in the proof of \textbf{Sub-case Two}, we have
		\begin{equation*}
			\bigg( \bigcup_{k_0\leq k \leq k_1 } C_k \bigg) \bigcap \{ \xi_2 | |\xi_2| \sim N_2  \} \bigcap S(r)
		\end{equation*}
		is contained in a strip which is parallel to $\R$-axis with width $N_1\theta$. Therefore, we obtain
		\begin{equation*}
			\big| \big( \tau_{k_1} - \frac{1}{2} \big) ^{\frac{1}{2}} - \big( \tau_{k_0} + \frac{1}{2} \big) ^{\frac{1}{2}} \big| \lesssim N_1 \theta,
		\end{equation*}
		which implies
		\begin{equation*}
			\tau_{k_1} - \tau_{k_0}  \lesssim \theta N_1 \tau_{k_1}^{\frac{1}{2}}.
		\end{equation*}
		When $\tau_{k_0} \le \frac{1}{100}(\tau_{k_1} - \tau_{k_0})$, we have
		\begin{equation*}
			\tau_{k_1} - \tau_{k_0} \lesssim \theta N_1 (\tau_{k_0} + \tau_{k_1} - \tau_{k_0})^{\frac{1}{2}} \lesssim \theta N_1 (\tau_{k_1} - \tau_{k_0})^{\frac{1}{2}},
		\end{equation*}
		which implies
		\begin{equation*}
			\tau_{k_1} - \tau_{k_0} \lesssim \theta^2 N_1^2.
		\end{equation*}
		Recall we have proved that $\{ \tau_{k} \}$ is an almost arithmetic sequence with common difference $N_1/\lambda$, we get
		\begin{equation*}
			k_1 - k_0 \lesssim  \lambda \theta^2 N_1.
		\end{equation*}
		Notice that $|C_k| \sim \lambda^{-1} \tau_k^{-\frac{1}{2}}$, we obtain
		\begin{align*}
			\sum_{k = k_0}^{k_1} |C_k| \lesssim & \sum_{k = k_0}^{k_1} \lambda^{-1}\tau_{k}^{-\frac{1}{2}} \lesssim \lambda^{-1} \sum_{k = k_0}^{k_1} \bigg( \tau_{k_0} + (k - k_0)\frac{N_1}{\lambda} \bigg)^{-\frac{1}{2}} \\
			\lesssim & \lambda^{-1} \tau_{k_0}^{-\frac{1}{2}} + \lambda^{-1} \sum_{k = k_0 + 1}^{k_1} \big[ (k-k_0)\frac{N_1}{\lambda} \big]^{-\frac{1}{2}}\\
			\lesssim & \frac{1}{\lambda M} + \theta.
		\end{align*}
		Here we have used the fact that $\lambda^{-1} \tau_{k_0}^{-\frac{1}{2}} \lesssim \frac{1}{\lambda M}$.
		When $\tau_{k_0} \ge \frac{1}{100} (\tau_{k_1} - \tau_{k_0})$, we get
		\begin{equation*}
			k_1 - k_0  \lesssim \lambda \theta \tau_{k_0}^{\frac{1}{2}}.
		\end{equation*}
		Then, we estimate as follows
		\begin{align*}
			\sum_{k = k_0}^{k_1} |C_k| \lesssim \sum_{k = k_0}^{k_1} \lambda^{-1}\tau_{k}^{-\frac{1}{2}} \lesssim (k_1 - k_0 + 1) \lambda^{-1}\tau_{k_0}^{-\frac{1}{2}}
			\lesssim \theta.
		\end{align*}
		
		We now turn to the \textbf{case $N_1 \sim N_2$}. Recall that the radius of annulus $S(r)$ can be approximated by
		\begin{equation}
			r \sim \frac{\sqrt{2\tau - |\xi|^2}}{2} \sim \frac{1}{2} \sqrt{2(2|\xi_2|^2 + |\xi|^2 - 2 |\xi||\xi_2| cos\angle(\xi_2 , \xi - \xi_2) ) }.
		\end{equation}
		By the angular condition $|cos\angle(\xi_2 , \xi - \xi_2)| \le \theta$, we have
		\begin{equation*}
			2|\xi_2|^2 + |\xi|^2 - 2 |\xi||\xi_2| cos\angle(\xi_2 , \xi - \xi_2) \ge |\xi_2|^2 \gtrsim N_1^2.
		\end{equation*}
		Therefore, we conclude that $r\sim N_1$. 
		
		As we have argued for the case when $N_2 \ll N_1$, we are aiming to show $K(\tau,\xi)$ is contained within portion of annulus $S(r)$ with arc length $O(\langle N_1\theta \rangle)$. In fact, assuming that $\angle(\xi_2 , \R) = \varphi + O(\theta)$ and denote
		\begin{align*}
			K^* := \sup \left\{ \eta \in \Z_{1/\lambda} : (\mu,\eta) \in K(\tau,\xi) \right\}\\
			K_* := \inf \left\{ \eta \in \Z_{1/\lambda} : (\mu,\eta) \in K(\tau,\xi) \right\}.
		\end{align*}
		On one hand, we can verify that
		\begin{equation*}
			K^* - K_* \lesssim \langle \lambda N_1 \theta \cos  \varphi \rangle.
		\end{equation*}
		On the other hand, by Figure 1 we have used in the previous case, for any $k \in \Z$
		\begin{equation*}
			\left| K(\tau,\xi) \cap \R\times\\{k/\lambda} \right| \lesssim \frac{1}{\lambda r \cos\varphi}.
		\end{equation*}
		Recall that $\frac{1}{r \cos\varphi} \le \frac{1}{ M}$, hence claim \eqref{fml-claim-coun} follows.
		
		Therefore, we have completed the estimate for \eqref{angular-reduce}. Now, we turn to the general bilinear form \eqref{reduce-2-angular} invoking the bilinear estimate for \eqref{angular-reduce}. 
		By angular truncation $| cos\angle (\xi_1 , \xi_2) |\le \theta$, we can find either $|arg(\xi_1) - arg(\xi_2)| = \pi/2 + O(\theta)$ or $|arg(\xi_1) - arg(\xi_2)| = 3\pi/2 + O(\theta)$. Hence, either $|l_1-l_2| = \pi/(2\theta) + O(1)$ or $|l_1-l_2| = 3\pi/(2\theta) + O(1)$. Using Minkowski's inequality, we obtain
		\begin{equation*}
			\|F\|_{L^2_{t,z}([0,1]\times\rt)} \lesssim \sum_{l_1,l_2: |l_1-l_2|= \pi/(2\theta) + O(1) \text{or} |l_1-l_2| = 3\pi/(2\theta) + O(1)} \bigg( \frac{1}{\lambda M} + \theta \bigg)^{\frac{1}{2}} \|\phi_{N_1,l_1}\|_{L^2_x} \|\phi_{N_2,l_2}\|_{L^2_x}.
		\end{equation*}
		Noticing that for fixed $l_1$, the above summation with respect to $l_2$ contains finite terms, thus by Cauchy-Schwarz's inequality and Plancherel's theorem we get
		\begin{align*}
			\|F\|_{L^2_{t,z}([0,1]\times\rt)} \lesssim & \bigg( \frac{1}{\lambda M} + \theta \bigg)^{\frac{1}{2}} \big( \sum_{l_1} \| \phi_{N_1,l_1} \|_{L^2_{x}}^2 \big)^{\frac{1}{2}}  \big( \sum_{l_2} \| \phi_{N_2,l_2} \|_{L^2_{x}}^2 \big)^{\frac{1}{2}}\\ 
			\lesssim & \bigg( \frac{1}{\lambda M} + \theta \bigg)^{\frac{1}{2}}  \|\phi_{N_1}\|_{L^2_{x}} \|\phi_{N_2}\|_{L^2_{x}}.
		\end{align*}
		Hence, we complete the proof.
	\end{proof}
	
	\subsubsection{More bilinear estimates to be used}	
	In the further application, we also need the refined angularly in the Bourgain space setting.
	
	\begin{proposition}\label{prop-ref-bil}
		Let $0 < \theta \ll 1 < M$ and $1 < N_2 \le N_1$. Suppose $u$ and $v \in X^{0,\frac{1}{2}+}$ are two functions which are frequency-localized around $N_1,N_2$ respectively. Then, the following bilinear form
		\begin{align*}
			F_1(t,z) := \int_{\R^2} \int_{(\rz)^2} & e^{it (\tau_1 + \tau_2) + iz \cdot (\xi_1 + \xi_2) } 1_{|\mu_1 - \mu_2|\gtrsim M}\nonumber \\
			&\times 1_{| cos\angle (\xi_1 , \xi_2) |\le \theta}\widehat{u}(\tau_1,\xi_1) \widehat{v}(\tau_2,\xi_2) (d\xi_1)_{\lambda} (d\xi_2)_{\lambda} d\tau_1 d\tau_2, 
		\end{align*}
		with $1\le M \lesssim N_1$ enjoys the bound
		\begin{equation*}
			\| F_1 \|_{L^2(\R \times \R \times \mathbb{T}_{\lambda})} \lesssim \left( \frac{1}{\lambda M} + \theta \right)^\frac12 \|u\|_{X^{0,\frac{1}{2}+}} \|v\|_{X^{0,\frac{1}{2}+}} .
		\end{equation*}
	\end{proposition}
	
	We have discussed refined bilinear estimate for very small $\theta$. However, we need the following estimate if the angular truncation is missing.
	\begin{proposition}\label{prop-no-theta}
		Let $1 < N_2 \le N_1$. Suppose $u$ and $v \in X^{0,\frac{1}{2}+}$ are two functions which are frequency-localized around $N_1,N_2$ respectively. Then, the following bilinear form
		\begin{align*}
			F_2(t,z) := \int_{\R^2} \int_{(\rz)^2} & e^{it (\tau_1 + \tau_2) + iz \cdot (\xi_1 + \xi_2) } 1_{|\mu_1 - \mu_2|\gtrsim M}\nonumber \\
			&\times \widehat{u}(\tau_1,\xi_1) \widehat{v}(\tau_2,\xi_2) (d\xi_1)_{\lambda} (d\xi_2)_{\lambda} d\tau_1 d\tau_2, 
		\end{align*}
		with $1\ll M \le N_1$ enjoys the bound
		\begin{equation*}
			\| F_2 \|_{L^2(\R \times \R \times \mathbb{T}_{\lambda})} \lesssim \bigg( \frac{1}{\lambda M} + \frac{N_2}{M} \bigg)^{\frac{1}{2}} \|u\|_{X^{0,\frac{1}{2}+}} \|v\|_{X^{0,\frac{1}{2}+}} .
		\end{equation*}
	\end{proposition}
	\begin{proof}
		Arguing as in Theorem \ref{bil-intro} and Proposition \ref{prop-ref-bil}, it is sufficient to prove
		\begin{equation*}
			\sup\limits_{(\tau , \xi)\in (\R \times \rz)}|M_1(\tau,\xi)| \lesssim \frac{N_2}{M},
		\end{equation*}
		where
		\begin{align}\label{fml-M1-unmbercont}
			M_1(\tau , \xi) := \{ & \xi_2 = (\mu_2 , \eta_2) \in \R \times \Z_{1/\lambda} : |\xi_2| \sim N_2, |\xi - \xi_2| \sim N_1, \nonumber \\ & | \mu/2 - \mu_2| \gtrsim M,  |\tau - |\xi_2|^2 - |\xi - \xi_2|^2 | \le 1  \}.
		\end{align}
		Similar with the proof in Theorem \ref{bil-intro}, we find that $M_1(\tau,\xi)$ is contained in an annulus $S(r)$, where
		\begin{equation*}
			S(r) = \bigg\{ \xi_2 : |\xi_2 - \frac{\xi}{2}| = r + O(1/r) \bigg\}, M\lesssim r\lesssim N_1.
		\end{equation*}
		We claim that at most $\lambda N_2$ discrete lines contained in set $M_1(\tau, \xi)$, since $|\xi_2| \sim N_2$. Moreover, $| \mu/2 - \mu_2| \gtrsim M$ shows that the measure of the intersection of each line with the annulus $S(r)$ is less than $\frac{1}{\lambda M}$. Hence, we have
		\begin{equation*}
			\sup\limits_{(\tau , \xi)\in (\R \times \rz)}|M_1(\tau,\xi)| \lesssim (1 + \lambda N_2) \frac{1}{\lambda M} = \frac{1}{\lambda M}+ \frac{N_2}{M}.
		\end{equation*}
		We have proved this proposition.
	\end{proof}
	\begin{remark}
		Proposition \ref{prop-no-theta} also holds if we replace $u$ with $\bar{u}$ and $|\mu_1 - \mu_2|$ with $|\mu_1 + \mu_2|$.
	\end{remark}

	In fact, Theorem \ref{bil-intro} gives refined bilinear estimate under two angular assumption, since $\xi_1 = \xi - \xi_2$ is almost determined by $\xi_2$. However, we also need a bilinear estimate with a weaker condition. Hence, we prove a bilinear estimate, which has only one angular restriction in the following.
	\begin{lemma}\label{lem-ref-bil-2}
		Let $0 < \theta \ll 1 < M$, $1 < N_2 \le N_1$ and $l\in \Z$. Suppose $\phi_{N_1}$ and $\phi_{N_2,l}$ are two functions whose frequencies are localized around $N_1$ and $N_2$. Moreover, the Fourier transform of $\phi_{N_2,l}$ is supported on the sector $\{ \xi \in \rz | |arg(\xi) - l\theta| \lesssim \theta \}$. Then, the following bilinear form
		\begin{align*}
			F_3(t,z) := \int_{(\rz)^2} & e^{-it (|\xi_1|^2 - |\xi_2|^2) - iz \cdot (\xi_1 + \xi_2) } 1_{|\mu_1 + \mu_2|\gtrsim M}  \\
			&\times \widehat{\phi_{N_1}}(\xi_1) \widehat{\phi}_{N_2,l}(\xi_2) (d\xi_1)_{\lambda} (d\xi_2)_{\lambda}, 
		\end{align*}
		enjoys the bound
		\begin{equation*}
			\| F_3 \|_{L^2([0,1] \times \R \times \mathbb{T}_{\lambda})} \lesssim \bigg( \frac{1}{\lambda M}+ \frac{N_1 \theta}{M}\bigg)^{\frac{1}{2}} \|\phi_{N_1}\|_{L^2} \|\phi_{N_2}\|_{L^2} .
		\end{equation*}
	\end{lemma}
	\begin{proof}
		Similar to the proof of Theorem \ref{bil-intro}, it is sufficient to show
		\begin{equation}\label{fml-tiltM-estimate}
			\sup\limits_{(\tau , \xi)\in (\R \times \rz)} |\widetilde{K}(\tau,\xi)| \lesssim \frac{1}{\lambda M} + \frac{N_1 \theta}{M} 
		\end{equation}
		where
		\begin{align*}
			\widetilde{K}(\tau , \xi) := \{ & \xi_2 = (\mu_2 , \eta_1) \in \R \times \Z_{1/\lambda} |\\ & |arg(\xi - \xi_2) - l\theta| \lesssim \theta, |\xi - \xi_2| \sim N_1, |\xi_2| \sim N_2 \\ 
			& |\mu| \sim M, |\tau - |\xi_2|^2 - |\xi - \xi_2|^2 | \le 1  \}.
		\end{align*}
		Compared to the set $M$ in Theorem \ref{bil-intro}, the set $\widetilde{M}$ here lacks an angular restriction. So the estimate we obtained here is rougher than Theorem \ref{bil-intro}. Fortunately, this estimate is sufficient to yield an improved energy increment. Similarly, we also divide $|\xi - \xi_2| = r \lesssim N_1$ and $|arg(\xi - \xi_2) - l\theta| \lesssim \theta$ into arcs of length $O(1/\lambda)$ with angular $O(1/(\lambda\langle r\rangle))$. Similarly, we have the number of the arcs can be bounded by $O(\langle \theta N_1 \rangle)$. Notice that 
		\begin{equation*}
			\big| \tau - |\xi_2|^2 + |\xi - \xi_2|^2  \big| = \big| \tau + |\xi|^2 - 2\mu \mu_2 - 2\eta \eta_2 \big| \lesssim 1, 
		\end{equation*}  
		and $|\xi| \gtrsim M$. Then, for fixed $\eta_2$, if $\xi_2,\xi_2^{\prime} \in K(\tau,\xi)$, we have $|\xi_2 - \xi_2^{\prime}| \lesssim \frac{1}{M}$. Therefore, we conclude that
		\begin{equation*}
			| \widetilde{K}(\tau,\xi) |\lesssim \frac{\langle \lambda N_1 \theta \rangle}{\lambda M} \lesssim \bigg( \frac{1}{\lambda M} + \frac{N_1 \theta}{M} \bigg).
		\end{equation*}
		We have finished the proof of this Lemma.
	\end{proof}
	
	Applying transfer principle to Lemma \ref{lem-ref-bil-2}, we arrive at
	\begin{proposition}\label{prop-ref-bil2}
		Let $0 < \theta \ll 1 < M$, $1 < N_2 \ll N_1$ and $l\in \Z$. Suppose $u,v \in X^{0,\frac{1}{2}+}$ are two functions which are frequency-localized around $N_1$ and $N_2$ respectively. Moreover, the Fourier transform of $v$ is supported on the sector $\{ \xi \in \rz | |arg(\xi) - l\theta| \lesssim \theta \}$. Then, the following bilinear form
		\begin{align*}
			F_4(t,z) := \int_{\R^2} \int_{(\rz)^2} & e^{it (\tau_1 + \tau_2) + iz \cdot (\xi_1 + \xi_2) } 1_{|\mu_1 - \mu_2|\gtrsim M}\\
			&\times \widehat{u}(\tau_1,\xi_1) \widehat{v}(\tau_2,\xi_2) (d\xi_1)_{\lambda} (d\xi_2)_{\lambda} d\tau_1 d\tau_2, 
		\end{align*}
		enjoys the bound
		\begin{equation*}
			\| F_4 \|_{L^2(\R \times \R \times \mathbb{T}_{\lambda})} \lesssim \bigg( \frac{1}{\lambda M} + \frac{N_1 \theta}{M} \bigg)^{\frac{1}{2}} \|u\|_{X^{0,\frac{1}{2}+}} \|v\|_{X^{0,\frac{1}{2}+}} .
		\end{equation*}
	\end{proposition}
	We also need the refined estimate for the product of linear propagator $e^{it\Delta}\phi_{N_2}\overline{e^{it\Delta}\phi_{N_2}}$.
	\begin{lemma}\label{lem-conjugate}
		Let $0 < \theta \ll 1 < M$, $1 < N_2 \ll N_1$ and $l\in \Z$. Suppose $u,v \in X^{0,\frac{1}{2}+}$ are two functions which are frequency-localized around $N_1$ and $N_2$ respectively. Moreover, the Fourier transform of $v$ is supported on the sector $\{ \xi \in \rz | |arg(\xi) - l\theta| \lesssim \theta \}$. Then, the following bilinear form
		\begin{align*}
			F_5(t,z) := \int_{\R^2} \int_{(\rz)^2} & e^{it (\tau_1 + \tau_2) + iz \cdot (\xi_1 + \xi_2) } 1_{|\mu_1 + \mu_2|\gtrsim M}\\
			&\times \widehat{u}(\tau_1,\xi_1) \widehat{\overline{v}}(\tau_2,\xi_2) (d\xi_1)_{\lambda} (d\xi_2)_{\lambda} d\tau_1 d\tau_2, 
		\end{align*}
		enjoys the bound
		\begin{equation*}
			\| F_5 \|_{L^2(\R \times \R \times \mathbb{T}_{\lambda})} \lesssim \bigg( \frac{1}{\lambda M} + \frac{N_1 \theta}{M} \bigg)^{\frac{1}{2}} \|u\|_{X^{0,\frac{1}{2}+}} \|v\|_{X^{0,\frac{1}{2}+}} .
		\end{equation*}
	\end{lemma}

\subsection{The proof of quadrilinear estimate}
	In this subsection, we prove the quadrilinear estimate by using the refined bilinear estimate established in the previous subsection.
	
	We define a multiplier as follows
	\begin{align*}
		A(\xi_1 , \xi_2 , \xi_3 , \xi_4) := & \big( m(\xi_1)\xi_1^2 - m(\xi_2)\xi_2^2 + m(\xi_3)\xi_3^2 - m(\xi_4)\xi_4^2 \big) \\
		& 1_{(\xi_1,\xi_2,\xi_3\xi_4) \in S^c}(\xi_1 , \xi_2 , \xi_3 , \xi_4),
	\end{align*}
	where $S^c$ is the complement of $S$ defined in $\eqref{fml-set-S}$. With this in hand, we  rewrite $\eqref{eq-6}$ as
	\begin{equation}\label{fml-rewrite-energy}
		(LHS)\mbox{ of }\eqref{eq-6}= \int_{0}^\delta \int_{\Sigma_4} A(\xi_1 , \xi_2 , \xi_3 , \xi_4) \hat{u}(t,\xi_1) \hat{\bar{u}}(t,\xi_2) \hat{u}(t,\xi_3) \hat{\bar{u}}(t,\xi_4) d\xi_1 d\xi_2 d\xi_3 d\xi_4 dt.
	\end{equation}
	From Lemma \ref{lemma-pointwise-multi}, multiplier $A$ have the following pointwise bound
	\begin{equation}\label{fml-pointbound-A}
		|A(\xi_1 , \xi_2 , \xi_3 , \xi_4)| \lesssim (min\{m_1,m_2,m_3,m_4\})^2 |\xi_{12}||\xi_{14}|.
	\end{equation}
	
	By using Littlewood-Paley decomposition, we rewrite
	\begin{equation*}
		\hat{u}(t,\xi_j) = \sum_{N_j \in 2^{\mathbb{N}}} \hat{u}_{N_j}(t,\xi_j), \quad j = 1,2,3,4.
	\end{equation*}
	Furthermore, we can assume $|\xi_j| \sim N_j\geq1$ with $j = 1,2,3,4$, and then we denote
	\begin{align*}
		J(N_1,N_2,N_3&,N_4) \\
		&:= \left|\int_{0}^\delta\int_{\Sigma_4} \Lambda_4(\xi_1,\xi_2,\xi_3,\xi_4)\widehat{u}_{N_1}(t,\xi_1)\widehat{\overline{u}}_{N_2}(t,\xi_2)\widehat{u}_{N_3}(t,\xi_3)\widehat{\overline{u}}_{N_4}(t,\xi_4)\,(d\xi_1)_\lambda\cdots\,(d\xi_4)_\lambda\,dt\right|.
	\end{align*}
	By the symmetry, we can allow $N_1 = max\{N_1,N_2,N_3,N_4\} \gtrsim N$ and $N_2 \ge N_4$. Hence, we only need to prove 
	\begin{equation*}
		|J(N_1,N_2,N_3,N_4)| \lesssim (N_1)^{0-} N^{-2+} \prod_{j=1}^4 \|Iu_{N_j}\|_{X^{1,\frac{1}{2}+}} .
	\end{equation*}
	In the following, we will consider two cases: $N_1 \sim N_3 \gg N_2 \gtrsim N_4$ and $N_1 \sim N_2 \gtrsim max \{N_3,N_4\}$. 
	
	\textbf{Case 1: $N_1 \sim N_3 \gg N_2 \gtrsim N_4$}. In this case we get $cos\angle (\xi_{12} , \xi_{14}) \sim 1 + O(N_2/N_1) \sim 1$. Recall that $A(\xi_1,\xi_2,\xi_3,\xi_4)$ is supported in the region $|cos\angle(\xi_{12},\xi_{14})| \ll 1$, therefore the multiplier vanishes.
	
	\textbf{Case 2: $N_1 \sim N_2 \gtrsim max \{N_3,N_4\}$}. By the symmetry, it is sufficient to consider when $N_1 \sim N_2 \gtrsim N_3 \ge N_4$. Therefore, Lemma \ref{lemma-pointwise-multi} yields
	\begin{equation}\label{fml-qua-pointwise}	|A(\xi_1,\xi_2,\xi_3,\xi_4)| \lesssim m(N_1)^2 N_3^2.
	\end{equation} 
	We further consider if $N_4$ is at least comparable with one of $N$ and $N_3$.
	
	\textbf{Case 2-1: $N_1 \sim N_2 \gtrsim N_3 \ge N_4, N_4\gtrsim min\{N,N_3\}$.} If $N_4 \gtrsim N$, We apply $L^4$ Strichartz estimate, $\eqref{fml-qua-pointwise}$, the fact $m(\xi)\lesssim1$ and get that
	\begin{align*}
		|J(N_1,N_2,N_3,N_4)| \lesssim & \frac{m(N_1)^2 N_3^2}{N_1 N_2 N_3 N_4} \prod_{j=1}^4 \| Iu_{N_j} \|_{X^{1,\frac{1}{2}+}} \\
		\lesssim & N_1^{0-} N^{-2+}.
	\end{align*}
	If $N_4 \sim N_3$, by the similar argument, we have
	\begin{align*}
		|J(N_1,N_2,N_3,N_4)| \lesssim & \frac{m(N_1)^2 N_3^2}{N_1 N_2 N_3 N_4} \prod_{j=1}^4 \| Iu_{N_j} \|_{X^{1,\frac{1}{2}+}} \\
		\lesssim & N_1^{0-} N^{-2+}.
	\end{align*}
	Then, we consider the case when $N_4$ is not comparable with $N$ nor $N_3$, that is $N_4 \ll \min\{N,N_3\}$.  
	
	\textbf{Case 2-2 : $N_1 \sim N_2 \gtrsim N_3$ and $N_4 \ll N_3$}. In this case, we can find
	\begin{equation*}
		| cos \angle (\xi_1 , \xi_3) | \lesssim | cos \angle (\xi_{14} , \xi_{34}) | + O(\frac{N_4}{N_3}) = N_1^{-1} + O(\frac{N_4}{N_3}).
	\end{equation*}
	Denote $\xi = (\mu , \eta) \in \R \times \Z_{1/\lambda}$, since $N_1 \ge N_3$ and $N_4 \ge 1$, one can find
	\begin{equation*}
		| cos \angle (\xi_1 , \xi_3) | = N_1^{-1} + O(\frac{N_4}{N_3}) = O\bigg(\frac{N_4}{N_3}\bigg) \sim \frac{\mu_1 \mu_3 + \eta_1  \eta_3}{N_1 N_3}.
	\end{equation*}
	Since $|cos \angle (\xi_1 , \xi_3)| \ll 1$, we have that $\xi_1$ is almost orthogonal with $\xi_3$. Therefore, at least one of the followings hold: $|\mu_1|\sim N_1$ and $|\mu_3| \sim N_3$. 
	
	Arguing as in \cite{Takaoka}, we further divide it into two subcases: $|\mu_1 + \mu_2| \ll N_3$ and $|\mu_1 + \mu_2| \gtrsim N_3$. If $|\mu_1 + \mu_2| \ll N_3$, we can find $|\mu_3| \ll N_3$ and $|\mu_1| \sim |\mu_2| \sim N_1$. Otherwise, if $|\mu_1 + \mu_2| \gtrsim N_3$, we have $|\mu_3| \sim N_3$. So, we will split it into two sub-cases $|\mu_1 + \mu_2| \ll N_3$ and $|\mu_1 + \mu_2| \gtrsim N_3$ in the following.
	
	\textbf{Case 2-2-1:  $N_1 \sim N_2 \gtrsim N_3$, $N_4 \ll N_3$, $|\mu_1 + \mu_2|\gtrsim N_3$ and $|\mu_3|\sim N_3$}. Using angular decomposition to function $u_{N_1}$, we obtain
	\begin{equation*}
		u_{N_1} = \sum_{l_1} u_{N_1,l_1}
	\end{equation*}
	and the Fourier transform of $u_{N_1,l_1}$ is supported in the sector $|arg(\xi_1) - l_1 \theta| \lesssim \theta$. Then, using angular decomposition $|arg(\xi_{34}) - l_{34}\theta| \lesssim \theta$ to  $u_{N_3}\bar{u}_{N_4}$, we have
	\begin{equation*}
		u_{N_3}\bar{u}_{N_4} = \sum_{l_{34}} F_{l_{34}},
	\end{equation*}
	where
	\begin{align*}
		F_{l_{34}} (t,z) = \int_{\R^2} \int_{(\rz)^2} e^{-it (\tau_3 + \tau_4) - iz \cdot (\xi_3 + \xi_4)} & 1_{|\mu_3| \sim N_3} 1_{arg(\xi_{34}) = l_{34}\theta + O(\theta)} \\
		&\hat{u}_{N_3}(\tau_3,\xi_3) \hat{\bar{u}}_{N_4}(\tau_4,\xi_4) (d\xi_3)_{\lambda} (d\xi_4)_{\lambda} d\tau_3 d\tau_4.
	\end{align*}
	Taking $\theta = O(\frac{N_4}{N_1})$, then we obtain
	\begin{equation*}
		| cos\angle (\xi_{12},\xi_1) - cos\angle (\xi_{12},\xi_{14}) | \lesssim O(\frac{N_4}{N_1}),
	\end{equation*}
	which implies 
	\begin{equation*}
		| cos \angle (\xi_{34} , \xi_1) | = |cos \angle (\xi_{12} , \xi_1)| \lesssim N_1^{-1} + \theta \lesssim \theta.
	\end{equation*}
	Therefore, $\xi_{34}$ is almost orthogonal with $\xi_1$ and we get either $|l_1 - l_{34}| = \frac{\pi}{2\theta} + O(1)$ or $|l_1 - l_{34}| = \frac{3\pi}{2\theta} + O(1)$ holds. In this case, by applying Proposition $\ref{prop-ref-bil}$ to $u_{N_1,l_1}\bar{u}_{N_2}$, Lemma $\ref{lem-conjugate}$ to $F_{l_{34}}$ and the fact that $\theta \sim \frac{N_4}{N_1}$, we arrive at
	\begin{align*}
		|J(N_1,N_2,N_3,N_4)| \lesssim & \sum_{\substack{|l_1 - l_{34}| = \pi/2\theta + O(1) \\ or |l_1 - l_{34}| = 3\pi/2\theta + O(1)}} m(N_1)^2N_3^2 \bigg( \frac{1}{\lambda N_3} + \theta \bigg)^{\frac{1}{2}} \\ 
		&\hspace{12ex} \times \|u_{N_1,l_1}\|_{X^{0,\frac{1}{2}+}}\|u_{N_2}\|_{X^{0,\frac{1}{2}+}} \|F_{l_{34}}\|_{L^2(\R\times\R\times\T_{\lambda})} \\
		\lesssim & \frac{m(N_1)^2N_3^2}{N_1 N_2 N_3 N_4} \bigg(\frac{1}{\lambda N_3} + \theta \bigg)^{\frac{1}{2}} \bigg(\frac{1}{\lambda N_3} + \frac{N_1 \theta}{N_3} \bigg)^{\frac{1}{2}} \\
		&\hspace{12ex} \times \|Iu_{N_1}\|_{X^{1,\frac{1}{2}+}}\|Iu_{N_2}\|_{X^{1,\frac{1}{2}+}} \|Iu_{N_3}\|_{X^{1,\frac{1}{2}+}} \|Iu_{N_4}\|_{X^{1,\frac{1}{2}+}} \\
		\lesssim & \bigg(\frac{N_3^{1/2}}{N_1^{3/2}N_2} + \frac{1}{\lambda N_1 N_2 N_4} 
		\bigg) \|Iu_{N_1}\|_{X^{1,\frac{1}{2}+}}\|Iu_{N_2}\|_{X^{1,\frac{1}{2}+}} \|Iu_{N_3}\|_{X^{1,\frac{1}{2}+}} \|Iu_{N_4}\|_{X^{1,\frac{1}{2}+}}\\
		\lesssim & N_1^{0-} N^{-2+} \|Iu_{N_1}\|_{X^{1,\frac{1}{2}+}}\|Iu_{N_2}\|_{X^{1,\frac{1}{2}+}} \|Iu_{N_3}\|_{X^{1,\frac{1}{2}+}} \|Iu_{N_4}\|_{X^{1,\frac{1}{2}+}}.
	\end{align*} 
	In this subcase, we obtain what we desired.

	\textbf{Case 2-2-2:  $N_1 \sim N_2 \gtrsim N_3$, $N_4 \ll N_3$, $|\mu_1 + \mu_2|\ll N_3$ and $|\mu_1| \sim N_1$}. In this case, $|\mu_1| \sim |\mu_2| \sim N_1$ and we obtain
	\begin{equation*}
		|\mu_1 + \mu_3| = |\mu_2 + \mu_4| \gtrsim N_1.
	\end{equation*}
	Taking $\theta \sim O(\frac{N_4}{N_3})$, we can verify that
	\begin{equation*}
		| cos \angle (\xi_{1} , \xi_3) | \le |cos \angle (\xi_{14} , \xi_{34})| + O(\frac{N_4}{N_3}) = |cos \angle (\xi_{14} , \xi_{12})| + O(\frac{N_4}{N_3}) \le N_1^{-1} + O(\frac{N_4}{N_3}) \sim O(\frac{N_4}{N_3}).
	\end{equation*}
	Moreover, we can find 
	\begin{equation*}
		|\mu_1 + \mu_3| = |\mu_2 + \mu_4| \sim N_2.
	\end{equation*}
	Applying Proposition $\ref{prop-no-theta}$ to $u_{N_1}u_{N_3}$ and Proposition \ref{prop-ref-bil} to $u_{N_2}u_{N_4}$ and $\theta \sim \frac{N_4}{N_3}$ we obtain
	\begin{align}\label{fml-qua-case2}
		|J(N_1,N_2,N_3,N_4)| \lesssim & m(N_1)^2N_3^2 \bigg( \frac{1}{\lambda N_1} + \theta \bigg)^{\frac{1}{2}} \bigg(\frac{1}{\lambda N_1} + \frac{N_4}{N_3} \bigg)^{\frac{1}{2}}  \prod_{j=1}^4 \|u_{N_j}\|_{X^{0,\frac{1}{2}+}} \\
		\lesssim & \frac{m(N_1)^2N_3^2}{N_1 N_2 N_3 N_4} \bigg( \frac{1}{\lambda N_1} + \theta \bigg)^{\frac{1}{2}} \bigg(\frac{1}{\lambda N_1} + \frac{N_4}{N_3} \bigg)^{\frac{1}{2}}  \prod_{j=1}^4 \|Iu_{N_j}\|_{X^{1,\frac{1}{2}+}} \\
		\lesssim &\bigg(\frac{1}{N_1 N_2} + \frac{1}{\lambda N_1 N_2 N_4}\bigg)    \prod_{j=1}^4 \|Iu_{N_j}\|_{X^{1,\frac{1}{2}+}} \\
		\lesssim & N_1^{0-} N^{-2+} \prod_{j=1}^4 \|Iu_{N_j}\|_{X^{1,\frac{1}{2}+}}.
	\end{align}
	So far, we have proved $\eqref{eq-6}$.

	\vspace{5mm}
	
\end{document}